\documentclass[preprint,1p,10pt]{elsarticle}

\journal{arXiv}

\biboptions{sort&compress,comma,round}

\usepackage{graphicx}
\usepackage{amssymb}
\usepackage{bm}
\usepackage[FIGTOPCAP]{subfigure}
\usepackage{amsmath}
\usepackage{float}

\DeclareSymbolFont{largesymbolsA}{U}{txexa}{m}{n}
\DeclareMathSymbol{\bigtimes}{\mathop}{largesymbolsA}{16}

\usepackage{color}
\usepackage[table]{xcolor}

\usepackage{setspace}

\usepackage{graphics}
\usepackage{amssymb, latexsym, amsmath, epsfig}
\usepackage{graphicx}
 \usepackage{color}
 \usepackage{epstopdf}
 \usepackage{comment}
 \usepackage{soul}
 \usepackage[normalem]{ulem}
  \usepackage{float}

\usepackage{amssymb,amsmath}
\usepackage{amsfonts,amstext,amsthm,amssymb,comment}
\usepackage{epsfig,graphics,color}

\graphicspath{{figures/}}

\numberwithin{equation}{section}

\newtheorem{remark}{Remark}

\newcommand{\bfs}[1]{{\boldsymbol #1}}

\newcommand{\Question}[1]{\marginpar{}}
\renewcommand{\Question}[1]{%
            \marginpar{\flushleft\scriptsize\bfseries\upshape#1}}

\begin{document}

\begin{frontmatter}


\title{A variationally separable splitting for the generalized-$\alpha$ method for parabolic equations}

%\author[add]{Ivo Babu\v{s}ka}
%

\author[ad1]{Pouria Behnoudfar}
\ead{pouria.behnoudfar@postgrad.curtin.edu.au}

\author[ad1,ad2]{Victor M. Calo}
\ead{Victor.Calo@curtin.edu.au}

\author[ad1,ad2]{Quanling Deng\corref{corr}}
\cortext[corr]{Corresponding author}
\ead{Quanling.Deng@curtin.edu.au}

\author[ad3]{Peter D. Minev}
\ead{minev@ualberta.ca}

%\address[ad]{Curtin Institute for Computation, Curtin University, Kent Street, Bentley, Perth, WA 6102, Australia}

\address[ad1]{Department of Applied Geology, Curtin University, Kent Street, Bentley, Perth, WA 6102, Australia}

\address[ad2]{Mineral Resources, Commonwealth Scientific and Industrial Research Organisation (CSIRO), Kensington, Perth, WA 6152, Australia}

\address[ad3]{Department of Mathematical and Statistical Sciences, University of Alberta, Edmonton, Alberta, Canada T6G 2G1}

\begin{abstract}
We present a variationally separable splitting technique for the generalized-$\alpha$ method for solving parabolic partial differential equations. We develop a technique for a tensor-product mesh which results in a solver with a linear cost with respect to the total number of degrees of freedom in the system for multi-dimensional problems. We consider finite elements and isogeometric analysis for the spatial discretization. The overall method maintains user-controlled high-frequency dissipation while minimizing unwanted low-frequency dissipation.  The method has second-order accuracy in time and optimal rates ($h^{p+1}$ in $L^2$ norm and $h^p$ in $L^2$ norm of $\nabla u$) in space. We present the spectrum analysis on the amplification matrix to establish that the method is unconditionally stable. Various numerical examples illustrate the performance of the overall methodology and show the optimal approximation accuracy.

\end{abstract}


\end{frontmatter}

\vspace{0.5cm}

\section{Introduction}

The generalized-$\alpha$ method was introduced by Chung and Hulbert in \cite{chung1993time} for solving (hyperbolic) structural dynamics problems. The method is of second-order accurate and possesses user-controlled numerical dissipation. The authors showed that the method renders high-frequency dissipation while minimizing unwanted dissipation low-frequency domain.  This method improves the Newmark methods \cite{newmark1959method} which possess numerical dissipation but are only first-order accurate and are too dissipative in the low-frequency region. The generalized-$\alpha$ method also improves the $\theta$ method of Wilson \cite{wilson1968computer}, the $\theta_1$ method of Hoff and Pahl \cite{hoff1988development}, and the $\rho$ method of Bazzi and Anderheggen \cite{bazzi1982rho}, which attain high-frequency dissipation with little low-frequency damping when maintaining second-order accuracy. Moreover, the method generalizes the HHT-$\alpha$ method of Hilber, Hughes, Taylor \cite{hilber1977improved} and the WBZ-$\alpha$ method of Wood, Bossak, and Zienkiewicz \cite{wood1980alpha}. That is, for particular values of the free parameters, the method reduces to HHT-$\alpha$ and WBZ-$\alpha$. The generalized-$\alpha$ method  produces an algorithm which provides an optimal combination of high-frequency and low-frequency dissipation in the sense that for a given value of high-frequency dissipation, the algorithm minimizes the low-frequency dissipation; see \cite{chung1993time}.

The generalized-$\alpha$ method was then extended to computational fluid dynamics governed by the parabolic-hyperbolic differential equations such as the Navier-Stokes equations in \cite{jansen2000generalized}. The method allows the user to control the high-frequency amplification factor. This parabolic time integrator encompasses a range of time integrators from zero damping midpoint rule to maximal damping Gear's method. Therein, the method was extended to filtered Navier-Stokes equations within the context of a stabilized finite element method.

One interesting feature of the time discretization of parabolic problems is that the self-adjoint and positive second-order spatial operator can be split in a sum of self-adjoint and positive operators of a simpler structure and
the original operator can be approximated via a factorized operator, that requires the solution of a set of much simpler linear problems.  Classical examples of such schemes are provided in the pioneering work of  Peaceman and Rachford \cite{peaceman1955numerical}, and Douglas and Rachford \cite{douglas1956numerical} that developed direction-splitting schemes for linear parabolic problems.   For a comprehensive review of this approach the reader is referred to the monographs of Marchuk  \cite{marchuk1990splitting} and Vabishchevich \cite{vabishchevich2013additive}.  Note that domain decomposition methods can also be interpreted in this framework as shown in \cite[Chapter 8]{SMV02}. In such a case, instead of splitting the operator
direction-wise as in the Peaceman-Rachford approach, it is decomposed into a sum of sub-operators by splitting the vector of unknowns using a partition of unity.  Although the resulting sub-operators have a different structure, the basic split formulation of 
these two approaches is very similar. Here we adopt the Peaceman-Rachford strategy, however, we adapt it to tensor-product approximations in the context of the finite element of isogeometric spatial approximations.
Note that similar splitting can be developed also for problems of a mixed hyperbolic-parabolic type, however, to the best of our knowledge, no systematic theory of stability and consistency has been developed in this case.

Splitting schemes developed on tensor-structured meshes reduce significantly the computational cost.  Ideally, the splitting solves the multi-dimensional problems with a computational time and storage which grow linearly  with respect to the number of degrees of freedom in the system. Based on tensor-product meshes, \cite{gao2014fast,gao2015preconditioners} apply alternating direction splitting schemes to reduce the computational cost of the resulting algebraic solver. The paper \cite{los2015dynamics} presents an application of alternating direction implicit splitting algorithm for
solving the parabolic equations using isogeometric finite element method. The authors show that
the overall scheme has a linear computational cost at every time step. The work \cite{los2017application} applies alternating direction splitting method using the isogeometric analysis to simulate tumor growth. Therein, the computational time cost scales linearly with the number of unknowns in the original multi-dimensional system.

We present a variationally separable splitting technique for the generalized-$\alpha$ method for parabolic equations. We formulate the variational formulation on tensor-product grids for multi-dimensional problems. Then the $d$-dimensional formulation is written as a product of $d$ formulations in each dimension plus error terms. We refer to these formulations as variationally separable. Based on the variational separability, we present a splitting technique for solving the resulting linear systems with a linear cost. With sufficient regularity, the approximate solution converges to the exact solution with optimal rates while reducing significantly the computational cost.

The rest of this paper is organized as follows. Section 2 describes the parabolic problem under consideration and introduces the particular spatial discretizations to arrive at the resulting matrix formulations of the problem. Section 3 presents a temporal discretization using the generalized-$\alpha$ method in an isogeometric analysis framework. Therein, we also introduce various splitting methods. Section 4 establishes the stability of the splitting schemes. We show numerically in Section 5 that the approximate solution converges optimally to the exact solution. We also verify that the computational cost is linear. Concluding remarks are given in Section 6.

\section{Problem statement and spatial discretization} \label{sec:ps}
Let $\Omega =(0, 1)^d \subset \mathbb{R}^d, d=2,3,$ be an open bounded domain. We consider the initial boundary-value problem for the heat equation
\begin{equation} \label{eq:pde}
\begin{aligned}
\frac{\partial u(x,t)}{ \partial t} - \Delta u(x, t) & = f(x,t), \qquad (x, t) \in \Omega \times (0, T], \\
u(x, t) & = u_D, \qquad \quad x \in \partial \Omega, \\
u(x, 0) & = u_0, \qquad \quad x \in \Omega, \\
\end{aligned}
\end{equation}
where the source $f$, the initial data $u_0$, and the Drichilet boundary data $u_D$ are assumed regular enough so that the problem admits a weak solution. For the sake of simplicity, we assume $u_D = 0$.

To simplify the arguments and given that we use tensor-product B-splines in multiple dimensions, we assume the partition $\mathcal{T}_h$ is a partition of $\Omega$ into non-overlapping tensor-product mesh elements.   Let $K\in \mathcal{T}_h$ be a generic element and its boundary is denoted as $\partial K$. Let $h$ denote the maximal diameter of the element $K$.  
%Let $\bfs{n}$ be the outward unit normal vector. 
Let $(\cdot, \cdot)_S$ denote the $L^2(S)$ inner product where $S$ is a $d$-dimensional domain ($S$ is typically $\Omega, K, \partial \Omega, \partial K$).

Next we define the discrete  space associated with the partition $\mathcal{T}_h$. For this purpose, we use the Cox-de Boor recursion formula~\cite{de1978practical,piegl2012nurbs} in each direction and then take tensor-product to obtain the necessary basis functions for multiple dimensions. If $X = \{x_0, x_1, \cdots, x_m \}$ is a knot vector with knots $x_j$,  then the $j$-th B-spline basis function of degree $p$, denoted as $\theta^p_j(x)$, is defined as~\cite{de1978practical, piegl2012nurbs}
\begin{equation} \label{eq:B-spline}
\begin{aligned}
\theta^0_j(x) & = 
\begin{cases}
1, \quad \text{if} \ x_j \le x < x_{j+1} \\
0, \quad \text{otherwise} \\
\end{cases} \\ 
\theta^p_j(x) & = \frac{x - x_j}{x_{j+p} - x_j} \theta^{p-1}_j(x) + \frac{x_{j+p+1} - x}{x_{j+p+1} - x_{j+1}} \theta_{j+1}^{p-1}(x).
\end{aligned}
\end{equation}

The Cox-de Boor recursion formula generates for a given knot vector a set of $C^k$ and $p$-th order B-spline basis functions, where $p=1,2,\cdots,$ and $k=0,1,\cdots, p-1$. The span of these basis functions generate a finite-dimensional subspace of the $H^1(\Omega)$ (see \cite{buffa2010isogeometric,evans2013isogeometric} for details):
\begin{equation} \label{eq:bs}
V^h = \text{span} \{ \Theta_j^p \}_{j=1}^{N_h} = 
\begin{cases}
S^p_k = \text{span} \{ \theta_j^p(x) \}_{j=1}^{N_x}, & \text{in 1D}\\
S^{p, q}_{k, m} = \text{span} \{ \theta_i^p(x) \theta_j^q(y) \}_{i, j=1}^{N_x, N_y}, & \text{in 2D}\\
S^{p, q,r}_{k, m,n} = \text{span} \{ \theta_i^p(x) \theta_j^q(y) \theta_l^r(z) \}_{i, j,l=1}^{N_x, N_y,N_z}, & \text{in 3D}\\
\end{cases}
\end{equation}
where $p,q,r$ and $k,m,n$ specifies the approximation order and continuity orders in each dimension, respectively. $N_x, N_y, N_z$ is the total number of basis functions in each dimension.

The weak formulation of the problem clearly reads: {\em Find $u \in H^1_0(\Omega)$ s.t.:}
\begin{equation} \label{eq:wf}
a(w, \dot u) + b(w, u) = \ell(w), \qquad \forall w \in H^1_0(\Omega), \ t>0,
\end{equation}
where $\dot u = \frac{\partial u}{\partial t}$ and 
\begin{equation}
a(w, v) = (w, v)_\Omega, \quad  b(w, v) =  (\nabla w, \nabla u)_\Omega, \quad \ell(w) = (w, f)_\Omega.
\end{equation}
The corresponding semi-discrete formulation is given by:  {\em Find $u_h(t) = u_h(\cdot, t) \in V^h$ s. t. for $t>0$:} 
\begin{equation} \label{eq:dwf}
a(w_h, \dot u_h) + b(w_h, u_h) = \ell(w_h), \qquad w_h \in V^h
\end{equation}
with $u_h(0)$ being the interpolant of $u_0$ in $V^h$. The discrete problem can be written in a matrix form as:
\begin{equation} \label{eq:mp}
M \dot{U} + KU = F,
\end{equation}
where $M$ and $K$ are the mass and stiffness matrices, $U$ is the vector of the unknowns, and $F$ is the source vector. The initial condition is given by
\begin{equation} \label{eq:u0}
U(0) = U_0,
\end{equation}
where $U_0$ is the given vector of initial condition $u_{0,h}$.  Since the space $V^h$ is constructed from tensor-product basis functions, the bilinear forms $a(\cdot, \cdot)$ and $b(\cdot, \cdot)$ inherit the tensor-product structure of the basis in some sense. For example, in 2D, we have that:
\begin{equation} \label{eq:sepa2d}
a(\theta_{i_x}^p(x) \theta_{i_y}^q(y), \theta_{j_x}^p(x) \theta_{j_y}^q(y)) = a_x(\theta_{i_x}^p(x), \theta_{j_x}^p(x) ) \cdot a_y(\theta_{i_y}^q(y),  \theta_{j_y}^q(y)),
\end{equation}
where:
\begin{equation}
a_x(\theta_{i_x}^p(x), \theta_{j_x}^p(x) ) = \int_0^1 \theta_{i_x}^p(x) \theta_{j_x}^p(x) \ \text{d} x, 
\quad a_y(\theta_{i_y}^q(y),  \theta_{j_y}^q(y)) = \int_0^1 \theta_{i_y}^q(y)  \theta_{j_y}^q(y) \ \text{d} y \\
\end{equation}
and:
\begin{equation} \label{eq:sepb2d}
\begin{aligned}
b(\theta_{i_x}^p(x) \theta_{i_y}^q(y), \theta_{j_x}^p(x) \theta_{j_y}^q(y)) & = b_x(\theta_{i_x}^p(x), \theta_{j_x}^p(x) ) \cdot a_y(\theta_{i_y}^q(y),  \theta_{j_y}^q(y)), \\
& \quad + a_x(\theta_{i_x}^p(x), \theta_{j_x}^p(x) ) \cdot b_y(\theta_{i_y}^q(y),  \theta_{j_y}^q(y)),
\end{aligned}
\end{equation}
where:
\begin{equation}
\begin{aligned}
b_x(\theta_{i_x}^p(x), \theta_{j_x}^p(x) ) & = \int_0^1 \frac{\text{d}}{\text{d} x} \theta_{i_x}^p(x) \frac{\text{d}}{\text{d} x} \theta_{j_x}^p(x) \ \text{d} x,  \\
\quad a_y(\theta_{i_y}^q(y),  \theta_{j_y}^q(y)) & = \int_0^1 \frac{\text{d}}{\text{d} y} \theta_{i_y}^q(y) \frac{\text{d}}{\text{d} y}  \theta_{j_y}^q(y) \ \text{d} y. \\
\end{aligned}
\end{equation}

Using this property of the discretization, called further on variational separability, we rewrite \eqref{eq:mp} as:
\begin{equation} \label{eq:sepmp2d}
(M_x \otimes M_y) \dot{U} + (K_x \otimes M_y + M_x \otimes K_y) U = F.
\end{equation}
Similarly, in 3D, we can rewrite \eqref{eq:mp} as:
\begin{equation} \label{eq:sepmp3d}
(M_x \otimes M_y \otimes M_z) \dot{U} + (K_x \otimes M_y \otimes M_z + M_x \otimes K_y \otimes M_z + M_x \otimes M_y \otimes K_z) U = F.
\end{equation}
Here, $M_\xi$ and $K_\xi$ with $\xi = x, y, z$ are one-dimensional mass and stiffness matrices, respectively. These matrices, as well as $M$ and $K$, are symmetric positive definite, and this property is crucial for obtaining stability of the generalized-$\alpha$ splitting scheme presented in the next section.

\section{Splitting schemes based on the generalized-$\alpha$ method}
\subsection{The generalized-$\alpha$ method}
Let us first recall the formulation of the the generalized-$\alpha$ method.
Consider a uniform partitioning of the time interval $[0,T]$ with a grid size $\tau$: $0 = t_0 < t_1 < \cdots < t_N = T$ and denote by $U_n, V_n$  the approximations to $U(t_n), \dot U(t_n)$, respectively. The time marching scheme of the generalized-$\alpha$ method is given by (see  \cite{jansen2000generalized}): 
\begin{equation} \label{eq:galpha}
\begin{aligned}
M V_{n+\alpha_m} + K U_{n+\alpha_f} & = F_{n+\alpha_f}, \\
U_{n+1} & = U_n + \tau V_n + \tau \gamma \delta(V_n), \\
U_0 & = U(0), \\
V_0 & = M^{-1} (F_0 - K U_0), \\
\end{aligned}
\end{equation}
where
\begin{equation} \label{eq:mf}
\begin{aligned}
F_{n+\alpha_f} & = F(t_{n+\alpha_f}), \\
W_{n+\alpha_g} & = W_n + \alpha_g \delta(W_n), \quad W = U, V, \quad g=m, f,\\
\delta(W_n) & = W_{n+1} - W_n. \\
\end{aligned}
\end{equation}
Substituting \eqref{eq:mf} into the first equation in \eqref{eq:galpha} we readily obtain:
\begin{equation} \label{eq:av}
\alpha_m A \delta(V_n) = F_{n+\alpha_f} - KU_n - (M + \tau \alpha_f K) V_n,
\end{equation}
where
\begin{equation}\label{eq:A}
A = M + \eta K \qquad \text{with} \quad \eta = \frac{\tau \gamma \alpha_f}{\alpha_m}.
\end{equation}

As shown by \cite{jansen2000generalized}, the scheme is formally second order accurate if:
\begin{equation} \label{eq:gm}
\gamma = \frac{1}{2} + \alpha_m - \alpha_f.
\end{equation}

The generalized-$\alpha$ method consists of two steps. One first solves \eqref{eq:av} for $\delta(V_n)$. Then apply to the second equation in \eqref{eq:galpha} to solve $U_{n+1}$, which is the solution at next time level.
Alternatively, supplementing \eqref{eq:av} with the second equation in \eqref{eq:galpha}, we arrive at a matrix formulation of the generalized-$\alpha$ method
\begin{equation} \label{eq:mp0}
\begin{aligned}
\begin{bmatrix}
U^{n+1} \\
\tau V^{n+1}\\
\end{bmatrix}
& =
\begin{bmatrix}
I - \frac{\tau \gamma}{\alpha_m} A^{-1} K & I - \frac{ \gamma}{\alpha_m} A^{-1} (M + \tau \alpha_f K) \\
- \frac{\tau}{\alpha_m} A^{-1} K & I - \frac{ 1}{\alpha_m} A^{-1} (M + \tau \alpha_f K) \\
\end{bmatrix}
\begin{bmatrix}
U^n \\
\tau V^n
\end{bmatrix} +
\begin{bmatrix}
\frac{\tau \gamma}{\alpha_m} A^{-1} F_{n+\alpha_f} \\
\frac{\tau}{\alpha_m}  A^{-1} F_{n+\alpha_f}
\end{bmatrix},
\end{aligned}
\end{equation}
where $I$ is an identity matrix which matches the dimension.

\subsection{First splitting scheme} \label{sec:ss1}
The main idea of the proposed splitting schemes is based on the following identity:
\begin{equation} \label{eq:a2a}
\begin{aligned}
A  = M + \eta K  = (M_x +  \eta K_x) \otimes (M_y +  \eta K_y) - \eta^2 K_x \otimes K_y,
\end{aligned}
\end{equation}
where the last term is clearly  of  order of $\tau^2$. We refer to the last term as the splitting error term.
This allows us to approximate $A$ in \eqref{eq:A} by
\begin{equation} \label{eq:A0}
\tilde A = (M_x +  \eta K_x) \otimes (M_y +  \eta K_y)
\end{equation}
up to a second order truncation error.
Substituting $\tilde A$ instead of $A$ in \eqref{eq:mp0} we arrive at the generalized-$\alpha$ splitting scheme:
\begin{equation} \label{eq:avs}
\alpha_m  (M_x +  \eta K_x) \otimes (M_y +  \eta K_y) \delta(V_n) = F_{n+\alpha_f} - KU_n - (M + \tau \alpha_f K) V_n.
\end{equation}

 Similarly, in 3D, we approximate $A$ in \eqref{eq:A} by
\begin{equation} \label{eq:A03d}
\tilde A = (M_x +  \eta K_x) \otimes (M_y +  \eta K_y) \otimes (M_z +  \eta K_z).
\end{equation}

\begin{remark}
The cost of the solution of the linear system with a matrix $M_\xi + \eta K_\xi$, where $\xi = x,y,z$ is linear with respect to the number of degrees of freedom (see also \cite{los2017application,los2015dynamics}). 
\end{remark}

\subsection{Second splitting scheme}
Similarly, we can also split the matrix on the right-hand side of \eqref{eq:av}. Splitting on both sides does not reduce the computational cost further. However, this second splitting delivers more accurate approximations in our numerical experiments. 

Firstly, we denote
\begin{equation}
B = M + \zeta K, \qquad \zeta = \tau \alpha_f.
\end{equation}
Now we approximate $B$ by the splitting idea, denoted as $\tilde B$. Thus, in 2D,
\begin{equation} \label{eq:B0}
\tilde B = (M_x + \zeta K_x ) \otimes (M_y + \zeta K_y ).
\end{equation}
The 3D system is split using a similar procedure to  \eqref{eq:A03d}.

Alternatively, we can split the matrices on both sides in a different way. We first write \eqref{eq:av} as follows
\begin{equation}\label{eq:mod}
  \alpha _m A \delta(V_n)=F_{n+\alpha_f}  - KU_n-\frac{\alpha_m}{\gamma} \big(A+\frac{\gamma -\alpha_m}{\alpha_m}M \big)V_n
\end{equation}
We then approximate $A$ by using $\tilde A$ defined as \eqref{eq:A0} in 2D and as \eqref{eq:A03d} in 3D. Similarly, this modification leads to a second-order accurate scheme in time. From our numerical experiments, splitting both the left-hand and right-hand sides of the equations deliver more accurate approximations.

\section{Spectral analysis} \label{sec:ea}
In this section, we perform the stability analysis to establish that the proposed splitting schemes are unconditionally stable. We start the analysis with the standard generalized-$\alpha$ method. Throughout this section, we set $F=0$ as it does not reduce the generality of the stability analysis. 

\subsection{The generalized-$\alpha$ method}
To study the stability, we spectrally decompose the matrix $K$ with respect to $M$  (see for example \cite{horn1990matrix}) to obtain
\begin{equation} \label{eq:sd}
K = MPDP^{-1},
\end{equation}
where $D$ is a diagonal matrix with entries to be the eigenvalues of the generalized eigenvalue problem
\begin{equation} \label{eq:eigKM}
K v = \lambda Mv
\end{equation}
and $P$ is a matrix with all the columns being the eigenvectors. We assume that all the eigenvalues are sorted in ascending order and are listed in $D$ and the $j$-th column of $P$ is associated with the eigenvalue $\lambda_j = D_{jj}$.

Note that $ I = P I P^{-1}. $  Using \eqref{eq:sd} and \eqref{eq:A}, we calculate
\begin{equation}
\begin{aligned}
A^{-1} & = (M + \eta K)^{-1} \\
& = (M + \eta MPDP^{-1} )^{-1} \\
& = \Big( M P ( I + \eta D ) P^{-1} \Big)^{-1} \\
& = P ( I + \eta D )^{-1} P^{-1} M^{-1}.
\end{aligned}
\end{equation}
Finally, we obtain:
\begin{equation}
\begin{aligned}
A^{-1} K & = \Big( P ( I + \eta D )^{-1} P^{-1} M^{-1} \Big) \Big( MPDP^{-1} \Big) = P ( I + \eta D )^{-1} D P^{-1}, \\
A^{-1} M & = \Big( P ( I + \eta D )^{-1} P^{-1} M^{-1} \Big) M = P ( I + \eta D )^{-1} P^{-1}. \\
\end{aligned}
\end{equation}

If we denote $E = ( I + \eta D )^{-1}$, we can rewrite the amplification matrix in \eqref{eq:mp0} as:
\begin{equation}
\begin{aligned}
\Xi & =
\begin{bmatrix}
I - \frac{\tau \gamma}{\alpha_m} A^{-1} K & I - \frac{ \gamma}{\alpha_m} A^{-1} (M + \tau \alpha_f K) \\
- \frac{\tau}{\alpha_m} A^{-1} K & I - \frac{ 1}{\alpha_m} A^{-1} (M + \tau \alpha_f K) \\
\end{bmatrix} \\
& =
\begin{bmatrix}
P & \bfs{0} \\
\bfs{0} & P \\
\end{bmatrix}
\begin{bmatrix}
I -  \frac{\tau \gamma}{\alpha_m} E D & I - \frac{ \gamma}{\alpha_m} E (I + \tau \alpha_f D)  \\
 -  \frac{\tau }{\alpha_m} E D & I - \frac{ 1}{\alpha_m} E (I + \tau \alpha_f D) \\
\end{bmatrix}
\begin{bmatrix}
P^{-1} & \bfs{0} \\
\bfs{0} & P^{-1} \\
\end{bmatrix}.
\end{aligned}
\end{equation}

Thus, we have
\begin{equation} \label{eq:mp1}
\begin{aligned}
\begin{bmatrix}
U^n \\
\tau V^n\\
\end{bmatrix}
& =
\begin{bmatrix}
P & \bfs{0} \\
\bfs{0} & P \\
\end{bmatrix}
\begin{bmatrix}
I -  \frac{\tau \gamma}{\alpha_m} E D & I - \frac{ \gamma}{\alpha_m} E (I + \tau \alpha_f D)  \\
 -  \frac{\tau }{\alpha_m} E D & I - \frac{ 1}{\alpha_m} E (I + \tau \alpha_f D) \\
\end{bmatrix}^n
\begin{bmatrix}
P^{-1} & \bfs{0} \\
\bfs{0} & P^{-1} \\
\end{bmatrix}
\begin{bmatrix}
U^0 \\
\tau V^0\\
\end{bmatrix}.
\end{aligned}
\end{equation}

Let us denote the matrix raised to power $n$ by $\tilde \Xi$. Clearly, the method would be unconditionally stable if the spectral radius of this matrix is bounded by one.
For the sake of completeness, we repeat  the analysis of \cite{jansen2000generalized}, and first consider the two limiting cases for $\tau$: $\tau\rightarrow 0$ and $\tau \rightarrow \infty$.    

In the limit $\tau\to 0$, since $D$ is diagonal, $\tau D\rightarrow 0$ and $E \to I$. The matrix $\tilde \Xi$ becomes an upper triangular matrix with eigenvalues:
\begin{equation}
\lambda_1 = 1 \qquad \text{and} \qquad \lambda_2 = 1 - \frac{1}{\alpha_m}.
\end{equation}
Both of them have the multiplicity equal to the dimension of the matrix $K$. This leads to the condition:
\begin{equation} \label{eq:am}
\alpha_m \ge \frac{1}{2}.
\end{equation}
Note that since the spectral radius of the matrix equals one, the method is stable over finite time intervals but not A-stable.

In the other limit of an infinite time step, the matrix $\tilde \Xi$ becomes a lower triangular matrix with eigenvalues:
\begin{equation} \label{eq:eig2}
\lambda_1 = 1 - \frac{1}{\alpha_f}, \qquad \text{and} \qquad \lambda_2 = 1 - \frac{1}{\gamma}.
\end{equation}
Both of them have multiplicity equal to the dimension of the matrix $K$. For second-order accuracy in time, $\gamma$ must satisfy \eqref{eq:gm}. This yields the condition:
\begin{equation}
\alpha_m \ge \alpha_f \ge \frac{1}{2},
\end{equation}
 Thus, the method is unconditionally stable in the two limiting cases, provided that its parameters satisfy the condition above. In the second limiting case
 it is also A-stable.

In order to control the high-frequency damping, Chung and Hulbert \cite{chung1993time} proposed to express the two parameters $\alpha_m$ and $\alpha_f$  in terms of the spectral radius $\rho_\infty$ corresponding to an infinite time step.
Using \eqref{eq:gm} and setting each eigenvalue in \eqref{eq:eig2} to be equal to  $-\rho_\infty$, we readily obtain:
\begin{equation} \label{eq:hfd}
\alpha_m = \frac{1}{2} \Big( \frac{3 - \rho_\infty}{1+\rho_\infty} \Big), \qquad \alpha_f = \frac{1}{1+\rho_\infty}.
\end{equation}
Condition \eqref{eq:gm} clearly infers that $\gamma = \alpha_f$.
This leads to a second-order accurate, unconditionally stable, one-parameter family of methods with a specified high frequency
damping.

In case of finite time steps the eigenvalues of the amplification matrix are the solutions of:
\begin{equation}\label{eq:char}
\begin{aligned}
0 & = \det( \tilde \Xi - \tilde \lambda I) \\
& = \det
\begin{bmatrix}
\tilde \Xi_{11} - \tilde\lambda I & \tilde \Xi_{12}  \\
\tilde \Xi_{21} & \tilde \Xi_{22} - \tilde\lambda I \\
\end{bmatrix} \\
& =
\det (\tilde \Xi_{11} - \tilde\lambda I) \cdot \det \Big( \tilde \Xi_{22} - \tilde\lambda I - \tilde \Xi_{21} (\tilde \Xi_{11} - \tilde\lambda I)^{-1} \tilde \Xi_{12} \Big).
\end{aligned}
\end{equation}

where
\begin{equation}
\tilde \Xi =
\begin{bmatrix}
\tilde \Xi_{11} & \tilde \Xi_{12}\\
\tilde \Xi_{21} & \tilde \Xi_{22}\\
\end{bmatrix}
=
\begin{bmatrix}
I -  \frac{\tau \gamma}{\alpha_m} E D & I - \frac{ \gamma}{\alpha_m} E (I + \tau \alpha_f D)  \\
 -  \frac{\tau }{\alpha_m} E D & I - \frac{ 1}{\alpha_m} E (I + \tau \alpha_f D) \\
\end{bmatrix}
\end{equation}
is a $2\times 2$ block matrix, with each block being a diagonal matrix. 

The first part of the eigenvalues of $\tilde \Xi $ are defined by the equation $\det (\tilde \Xi_{11} - \tilde\lambda I) =0$,
where
 $\tilde \Xi_{11} = I -  \frac{\tau \gamma}{\alpha_m} E D$ is a diagonal matrix with diagonal entries given by:
\begin{equation}
1 - \frac{\tau \lambda_k \gamma}{ \alpha_m + \tau \lambda_k \gamma \alpha_f},
\end{equation}
with $\lambda_k$ being the $k-$th diagonal entry of $D$.
To guarantee stability we need that the absolute value of each of the eigenvalues of $\tilde \Xi $ is bounded by one and therefore:
\begin{equation}
-1 \le 1 - \frac{\tau \lambda_k \gamma}{ \alpha_m + \tau \lambda_k \gamma \alpha_f} \le 1,
\end{equation}
or
\begin{equation}
0 \le \frac{\tau \lambda_k\gamma}{ \alpha_m + \tau \lambda_k \gamma \alpha_f} \le 2.
\end{equation}
The left-hand-side inequality is satisfied since $D$ is a positive definite matrix and all parameters are non-negative. The right-side inequality can be rewritten as:
\begin{equation}
\tau \lambda_k \gamma (1 - 2 \alpha_f) \le 2 \alpha_m.
\end{equation}
Since $\tau \lambda \ge 0$, the condition $1 - 2 \alpha_f \le 0$ is sufficient to guarantee that the absolute value of this part of the spectrum of $\tilde \Xi $ is bounded by one.

Since each matrix in the second condition in \eqref{eq:char}
\begin{equation}
\det \Big( \tilde \Xi_{22} - \tilde\lambda_k I - \tilde \Xi_{21} (\tilde \Xi_{11} - \tilde\lambda_k I)^{-1} \tilde \Xi_{12} \Big) = 0.
\end{equation}
is a diagonal matrix, the rest of the spectrum of $\tilde \Xi $ is a solution of:
\begin{equation} \label{eq:lam0}
\begin{aligned}
0 & = \Big( 1 - \frac{\tau \lambda_k \gamma}{ \alpha_m + \tau \lambda_k \gamma \alpha_f} - \tilde \lambda_k \Big) \cdot \Big( 1 - \frac{1+ \tau \lambda_k \alpha_f}{ \alpha_m + \tau \lambda_k \gamma \alpha_f} - \tilde \lambda_k  \Big) \\
& \quad + \Big( \frac{\tau \lambda_k}{ \alpha_m + \tau \lambda_k \gamma \alpha_f} \Big) \cdot \Big( 1 - \frac{\gamma + \tau \lambda_k \gamma \alpha_f}{ \alpha_m + \tau \lambda_k \gamma \alpha_f} \Big).
\end{aligned}
\end{equation}
It is not difficult to verify that a sufficient condition, that guarantees that the absolute value of these eigenvalues is bounded by one, is given by:
\begin{equation} \label{eq:condmf}
\alpha_m \ge \alpha_f \ge \frac{1}{2}.
\end{equation}

\subsection{Stability of the splitting schemes}
Now we perform a stability analysis for the splitting scheme proposed in Section \ref{sec:ss1}. Similarly, we apply the spectral decomposition of one of the directional matrices $K_\xi$ with respect to its directional $M_\xi$ and arrive at
\begin{equation} \label{eq:sd0}
K_\xi = M_\xi P_\xi D_\xi P_\xi^{-1},
\end{equation}
where $D_\xi$ is a diagonal matrix with entries being the eigenvalues of the generalized eigenvalue problem
\begin{equation} \label{eq:eigKM0}
K_\xi v_\xi = \lambda_\xi M_\xi v_\xi
\end{equation}
and $P_\xi$ is a matrix with all the columns being the eigenvectors. Herein, $\xi = x,y,z$ specifies the coordinate directions. We assume that all the eigenvalues are sorted in ascending order and are listed in $D_\xi$ and the $j$-th column of $P_\xi$ is associated with the eigenvalue $\lambda_{\xi,j} = D_{\xi,jj}$. We perform the analysis for 2D splitting as follows.

Using \eqref{eq:sd0} and \eqref{eq:A0}, we now calculate:
\begin{equation}
\begin{aligned}
\tilde A^{-1} & = (M_x +  \eta K_x)^{-1} \otimes (M_y +  \eta K_y)^{-1} \\
& = (M_x + \eta M_x P_x D_x P_x^{-1} )^{-1} \otimes (M_y + \eta M_y P_y D_y P_y^{-1} )^{-1} \\
& = P_x E_x P_x^{-1} M_x^{-1} \otimes P_y E_y P_y^{-1} M_y^{-1},
\end{aligned}
\end{equation}
where:
\begin{equation}
E_\xi = ( I + \eta D_\xi )^{-1}, \qquad \xi = x, y.
\end{equation}

Similarly to the case of the unsplit scheme we have:
\begin{equation}\label{eq:multi}
\begin{aligned}
\tilde A^{-1} M & = \Big( P_x E_x  P_x^{-1} M_x^{-1} \otimes P_y E_y P_y^{-1} M_y^{-1} \Big) \cdot ( M_x \otimes M_y) \\
& = P_x E_x  P_x^{-1} \otimes P_y E_y P_y^{-1} \\
& = \Big( P_x \otimes P_y \Big) \cdot \Big( E_x \otimes E_y \Big) \cdot \Big( P_x^{-1} \otimes P_y^{-1} \Big), \\
\tilde A^{-1} K & = \Big( P_x E_x  P_x^{-1} M_x^{-1} \otimes P_y E_y P_y^{-1} M_y^{-1} \Big)  \\
& \qquad  \cdot \Big( M_x P_x D_x P_x^{-1} \otimes M_y + M_x \otimes M_y P_y D_y P_y^{-1} \Big) \\
& = P_x E_x  D_x P_x^{-1} \otimes P_y E_y P_y^{-1} + P_x E_x P_x^{-1} \otimes P_y E_y D_y P_y^{-1}  \\
& = \Big( P_x \otimes P_y \Big) \cdot \Big( E_x D_x \otimes E_y + E_x \otimes E_y D_y \Big) \cdot \Big( P_x^{-1} \otimes P_y^{-1} \Big). \\
\end{aligned}
\end{equation}

If we use the following notation:
\begin{equation}
\begin{aligned}
I & = P_x I_x P_x^{-1} \otimes P_y I_y P_y^{-1}  = \Big( P_x \otimes P_y \Big) \cdot \Big( I_x \otimes I_y \Big) \cdot \Big( P_x^{-1} \otimes P_y^{-1} \Big), \\
G & =  -\frac{\tau}{\alpha_m} \big( E_x D_x \otimes E_y + E_x \otimes E_y D_y \big).
\end{aligned}
\end{equation}
we can write the blocks of the  amplification matrix of the scheme in this case as follows:
\begin{equation}
\begin{aligned}
\Xi_{11} & = I - \frac{\tau \gamma}{\alpha_m} \tilde A^{-1} K = \Big( P_x \otimes P_y \Big) \cdot \Big( I_x \otimes I_y + \gamma G \Big) \cdot \Big( P_x^{-1} \otimes P_y^{-1} \Big), \\
\Xi_{21} &  = - \frac{\tau}{\alpha_m} \tilde A^{-1} K 
 = \Big( P_x \otimes P_y \Big) \cdot G \cdot \Big( P_x^{-1} \otimes P_y^{-1} \Big) \\
\Xi_{12} & =
I - \frac{ \gamma}{\alpha_m} \tilde A^{-1} (M + \tau \alpha_f K) \\
& = \Big( P_x \otimes P_y \Big) \cdot \Big( I_x \otimes I_y - \frac{ \gamma }{\alpha_m} E_x \otimes E_y +\gamma \alpha_f G \Big) \cdot \Big( P_x^{-1} \otimes P_y^{-1} \Big), \\
\Xi_{22} & = I - \frac{ 1}{\alpha_m} \tilde A^{-1} (M + \tau \alpha_f K) \\
& = \Big( P_x \otimes P_y \Big) \cdot \Big( I_x \otimes I_y - \frac{ 1 }{\alpha_m} E_x \otimes E_y + \alpha_f G \Big) \cdot \Big( P_x^{-1} \otimes P_y^{-1} \Big). \\
\end{aligned}
\end{equation}
and the matrix itself as:
\begin{equation}\label{eq:stability}
\begin{aligned}
\Xi & =
\begin{bmatrix}
I - \frac{\tau \gamma}{\alpha_m} \tilde A^{-1} K & I - \frac{ \gamma}{\alpha_m} \tilde A^{-1} (M + \tau \alpha_f K) \\
- \frac{\tau}{\alpha_m} \tilde A^{-1} K & I - \frac{ 1}{\alpha_m} \tilde A^{-1} (M + \tau \alpha_f K) \\
\end{bmatrix} \\
& =
\begin{bmatrix}
P_x \otimes P_y & \bfs{0} \\
\bfs{0} & P_x \otimes P_y \\
\end{bmatrix}
\begin{bmatrix}
I_x \otimes I_y + \gamma G & I_x \otimes I_y - \frac{ \gamma }{\alpha_m} E_x \otimes E_y +\gamma \alpha_f G  \\
G & I_x \otimes I_y - \frac{ 1 }{\alpha_m} E_x \otimes E_y + \alpha_f G \\
\end{bmatrix} \\
& \qquad
\begin{bmatrix}
P_x^{-1} \otimes P_y^{-1} & \bfs{0} \\
\bfs{0} & P_x^{-1} \otimes P_y^{-1} \\
\end{bmatrix}.
\end{aligned}
\end{equation}
Obviously, the stability of the scheme is determined by the spectral radius of:
\begin{equation}
\tilde{\Xi}  =
\begin{bmatrix}
I_x \otimes I_y + \gamma G & I_x \otimes I_y - \frac{ \gamma }{\alpha_m} E_x \otimes E_y +\gamma \alpha_f G  \\
G & I_x \otimes I_y - \frac{ 1 }{\alpha_m} E_x \otimes E_y + \alpha_f G \\
\end{bmatrix}.
\end{equation}
In the limit $\tau \to 0$, we arrive at the same sufficient condition as in \eqref{eq:am}.
In the limit $\tau \to \infty$, we obtain the eigenvalues
$
\lambda_1 = 1,
$
with a multiplicity equal to twice the dimension of the matrix $K$.
This implies that in the limiting case $\tau = \infty$, the method is stable but not A-stable. Following similar arguments as in the case of the unsplit generalized-$\alpha$ method, one can show that the scheme is stable for any finite time step size. The stability analysis for splitting on both sides follow similar arguments and we omit the details herein for simplicity.

\begin{remark}
The stability analysis for 3D splitting is more involved but follows the same logic. We omit the details herein and state that the conditions on $\alpha_m, \alpha_f$ for unconditional stability are the same as for the 2D splitting.
\end{remark}

Lastly, before closing this section, we present error estimates. The normal modes analysis above establishes the stability of the splitting schemes. 
Since the splitting error, is formally second-order accurate in time, we may expect
that the splitting schemes are second-order accurate in time, provided that the exact solution is sufficiently regular (note that the stability of the splitting
requires a much higher degree of regularity of the exact solution). 
Thus, we expect to obtain estimate similar to the classical results for parabolic problems (see for example \cite{thomee1984galerkin}):
\begin{equation} \label{eq:errl2h1}
\begin{aligned}
\| u_h^n - u(t_n) \|_{0, \Omega} & \le C(u) (h^{p+1} + \tau^2), \\
\| \nabla u(T) - \nabla u_h(T)\|_{0,\Omega} & \le C(u) (h^p + \tau^2), \\
\end{aligned}
\end{equation}
where $u_h^n$ is the approximate solution at time $t_n$ and $C(u)$ is a positive constant independent of mesh size $h$ and time step $\tau$, and $p$ is the order of the spatial approximation.

\section{Numerical experiments} \label{sec:num}
In this section, we present numerical examples to show the performance of the proposed splitting schemes.  We focus on 2D splitting but also show some 3D results. The goal of the numerical results
presented below is to validate that the schemes result in optimal convergence rates in space and time,  and that the computational cost of the splitting schemes is linear with respect to the total number of degrees of freedom in the system. For this purpose, we consider \eqref{eq:pde} with an exact manufactured solution:
\begin{equation}\label{pro}
u = 
\begin{cases}
 u(x,y,t)=\sin(\pi x)\sin(\pi y) e^{-2\pi ^2 t}, \quad \text{in} \quad 2D, \\
 u(x,y,z,t)=\sin(\pi x)\sin(\pi y) \sin(\pi z)  e^{-3\pi ^2 t}, \quad \text{in} \quad 3D, \\
 \end{cases}
\end{equation}
from which the forcing function and boundary and initial conditions are derived.

\begin{figure}[ht]
\centering\includegraphics[width=6.5cm]{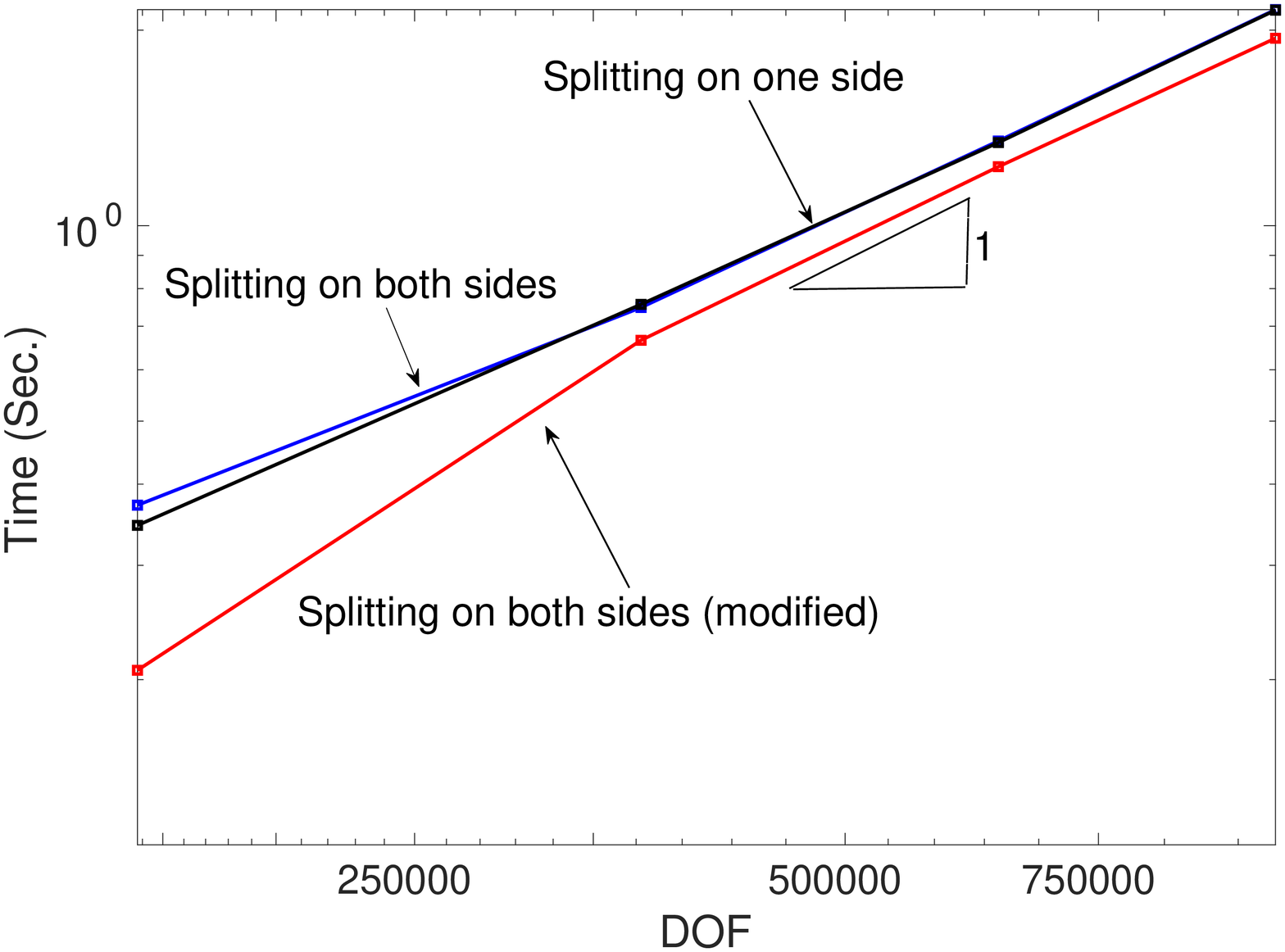} 
\centering\includegraphics[width=6.5cm]{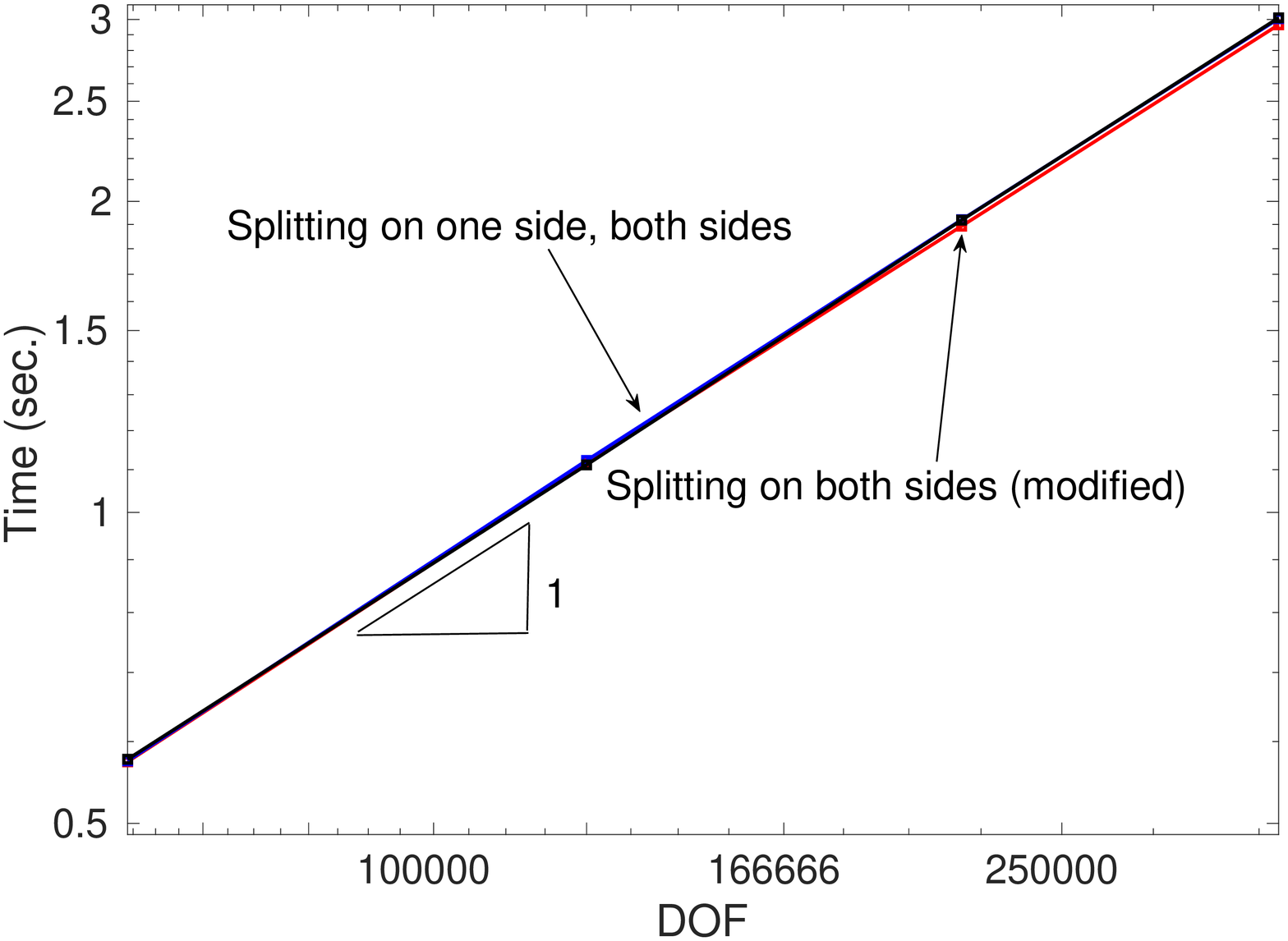} 
\caption{The computational cost (linear) of the proposed splitting schemes when using $C^1$ quadratic isogeometric elements with $\tau=10^{-3}$ and $\rho_\infty=0$ in 2D (left) and in 3D (right).}
 \label{fig:cost}
\end{figure}

First, we consider the computational cost for the parabolic problem \eqref{eq:pde} in both 2D and 3D.

For spatial discretizations with $p$-th order finite elements or isogeometric elements, the one-dimensional mass matrix $M_\xi$ and stiffness matrix $K_\xi, \xi = x, y, z,$ are of a half-bandwidth $p$ and so is $M_\xi + \eta K_\xi$. Assuming $K_\xi$ has a dimension $m_\xi$, then solving a linear matrix system with $M_\xi + \eta K_\xi$ being the matrix requires $\mathcal{O}(p^2m_\xi)$ operations when using Gaussian elimination. The main cost for solving \eqref{eq:av} or equivalently \eqref{eq:mp0} with the splitting schemes is on the inversion of the tensor-product structure matrix $\tilde A$. This cost requires $\mathcal{O}(p^2m_x m_y)$ operations for 2D problem and $\mathcal{O}(p^2m_x m_y m_z)$ operations for 3D problem. 
Thus, the cost for solving \eqref{eq:av} using splitting schemes grow linearly with respect to the degrees of freedom. The property remains for higher dimensional problems. Additionally,
it allows for the use of direct solvers for problems of any dimension. The reader is referred to \cite{los2017application} for more details on the linear solver.

Figure \ref{fig:cost} shows that the costs for solving the resulting algebraic matrix problems of all the proposed splitting schemes are linear with respect to the total number of degrees of freedom in the system for multi-dimensional problems.  Herein, we use a direct solver (Gaussian elimination) and as an example, we use $C^1$ quadratic isogeometric elements for the spatial discretization and a time step size $10^{-3}$, however,  we observe the same linear cost when using finite elements and isogeometric elements. This validates the efficiency of the splitting schemes when solving the resulting matrix problems. 

\begin{figure}[ht]
\centering\includegraphics[width=6.5cm]{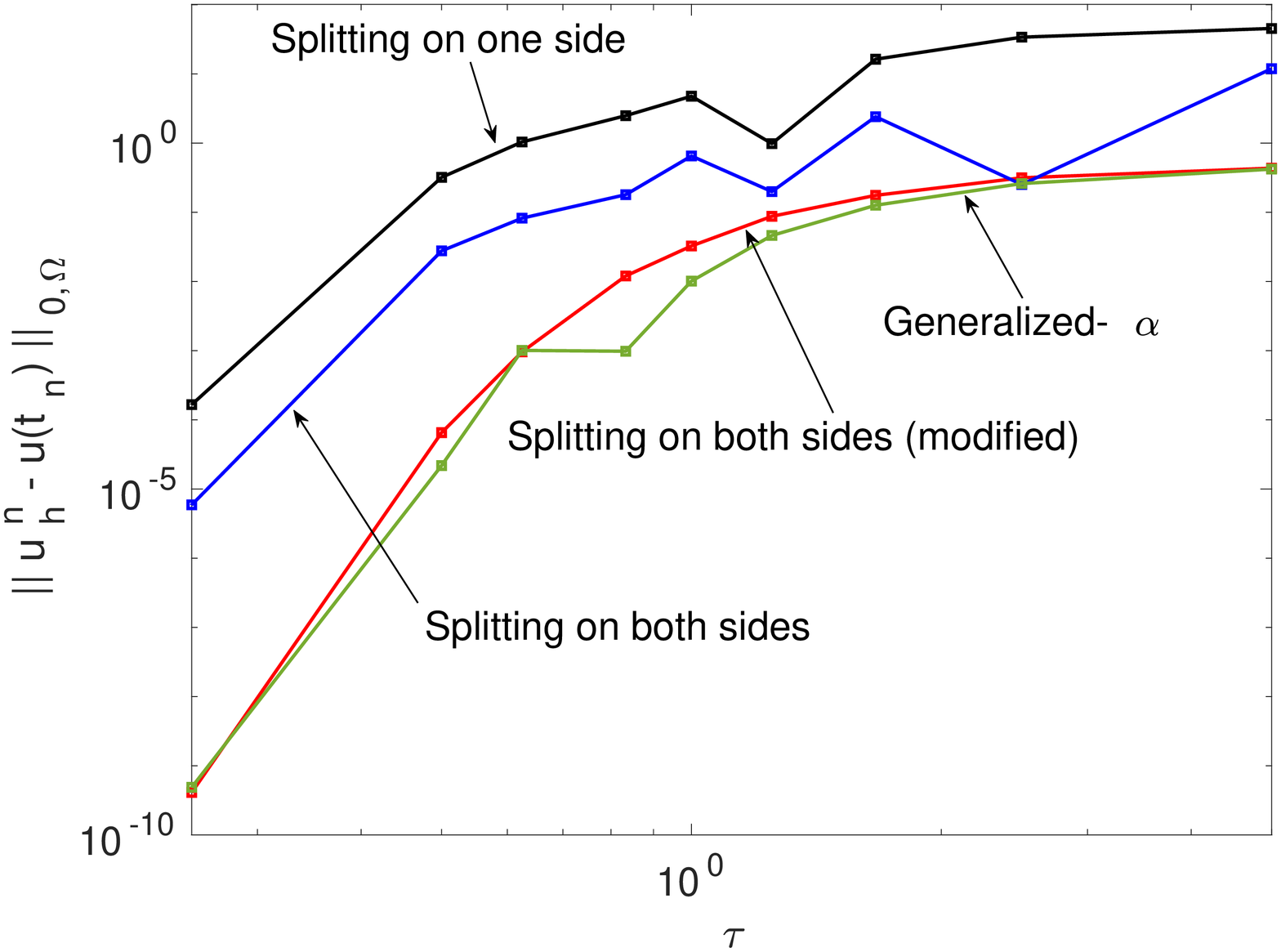} 
\centering\includegraphics[width=6.5cm]{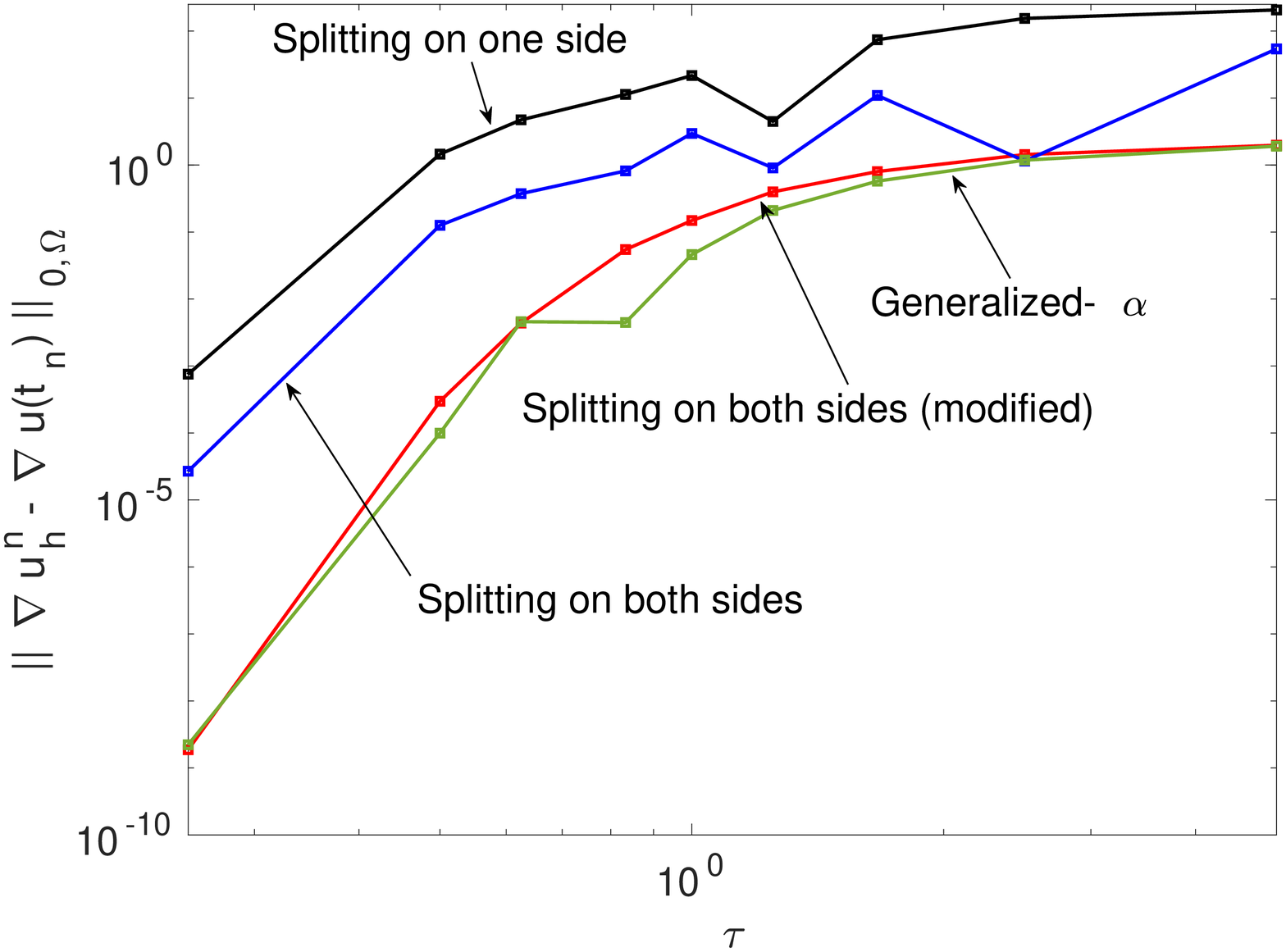} 
\caption{Stability validation on the proposed splitting schemes when using $C^1$ quadratic isogeometric elements with final time $T=5$, $\rho_\infty=0.5$, and $64 \times 64$ uniform elements in 2D.}
	\label{fig:stab}
\end{figure}

In Figure \ref{fig:stab} we present the $L^2$ norm of $u$ and $\nabla u$ errors at the final time $T=5$ with respect to time step size $\tau$. As $\tau$ increases, both errors approach a finite number, which validates numerically the unconditional stability of the generalized-$\alpha$ and splitting schemes. 
The scheme named "Splitting on both sides" refers to the scheme given by \eqref{eq:B0} while the "Splitting on both sides (modified)" refers to the scheme in equation \eqref{eq:mod}.
We observe the same behaviour for all the schemes with other scenarios, such as different $\rho_\infty$, mesh configurations, and finite elements with higher order basis functions. 

\begin{figure}[ht]
\subfigure[FEM $(p=2, C^0)$]{\centering\includegraphics[width=6.5cm]{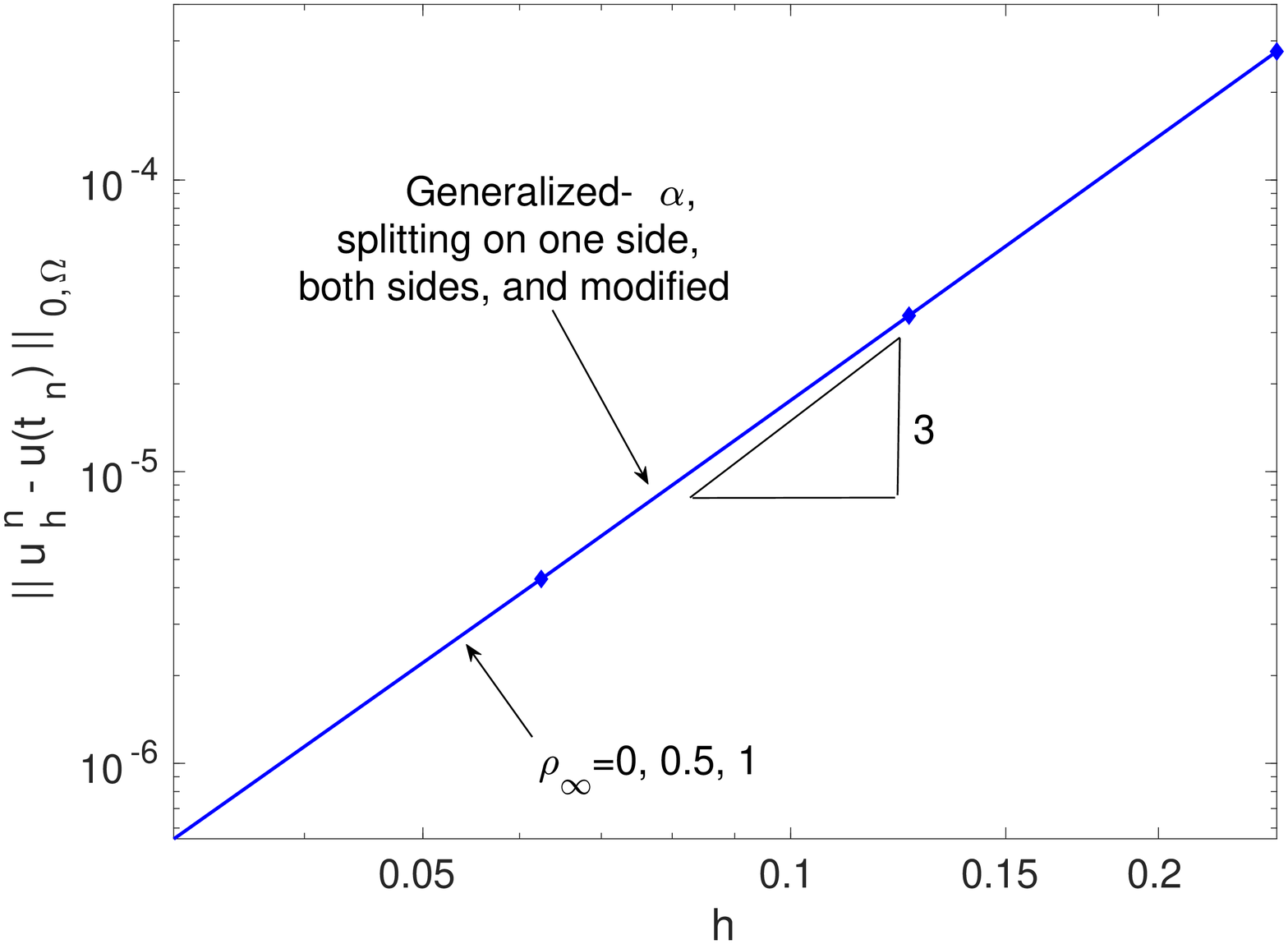}}
\subfigure[IGA $(p=2, C^1)$]{\centering\includegraphics[width=6.5cm]{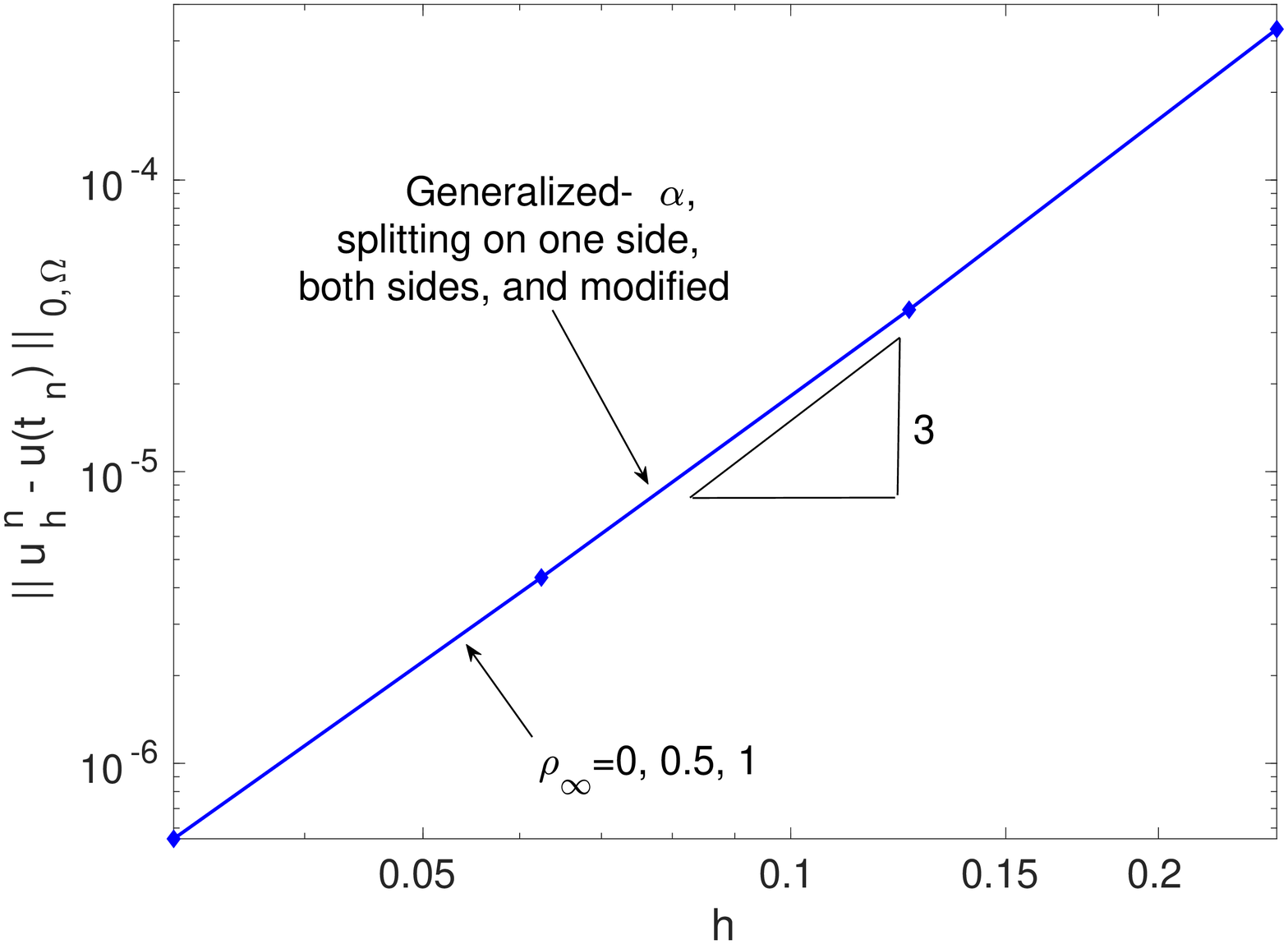}}
\subfigure[FEM $(p=2, C^0)$]{\centering\includegraphics[width=6.5cm]{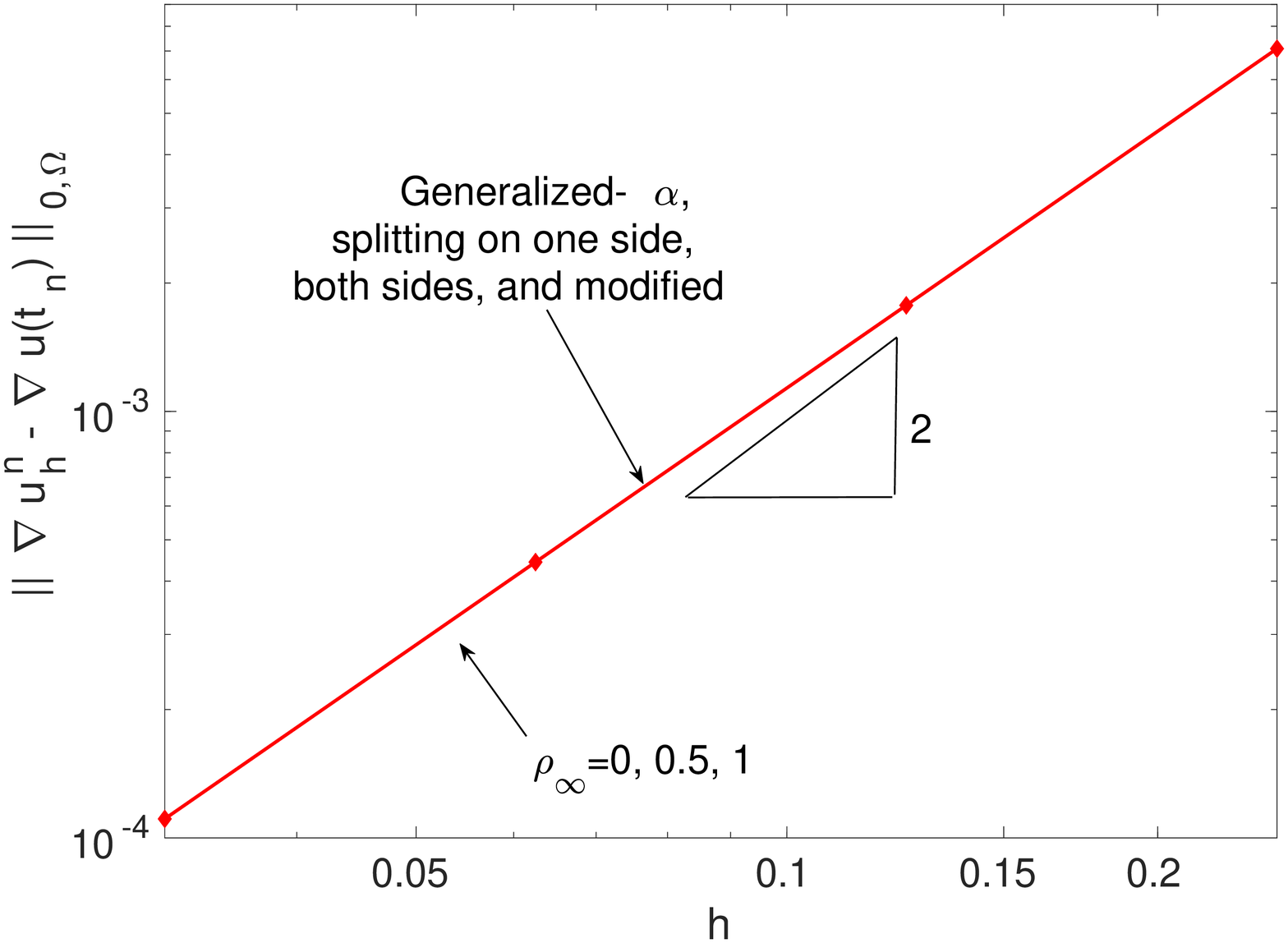}}
\hspace{0.3 cm}
\subfigure[IGA $(p=2, C^1)$]{\centering\includegraphics[width=6.5cm]{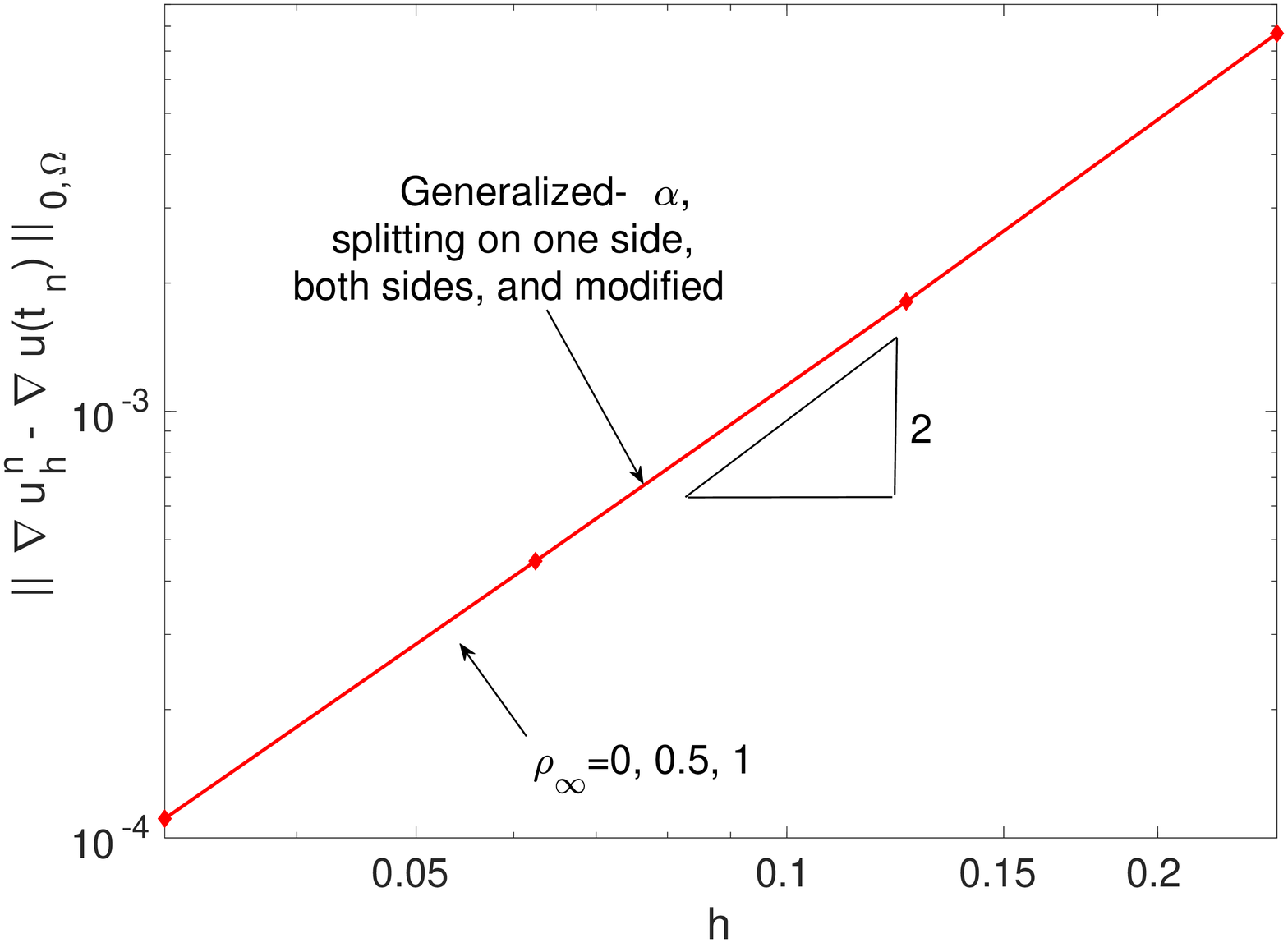}}
\caption{$L^2 $ norm of $u$ and $\nabla u$ errors when using classical $C^0$ quadratic finite elements and $C^1$ quadratic isogeometric elements for space discretization with $\rho_\infty=0, 0.5, 1$. The final time $t_n=T=1$ and time step-size is $\tau = 10^{-4}$.}
\label{fig:p2h}
\end{figure}

\begin{figure}[!ht]
\hspace{-0.12 cm}
\centering\includegraphics[width=6.5cm]{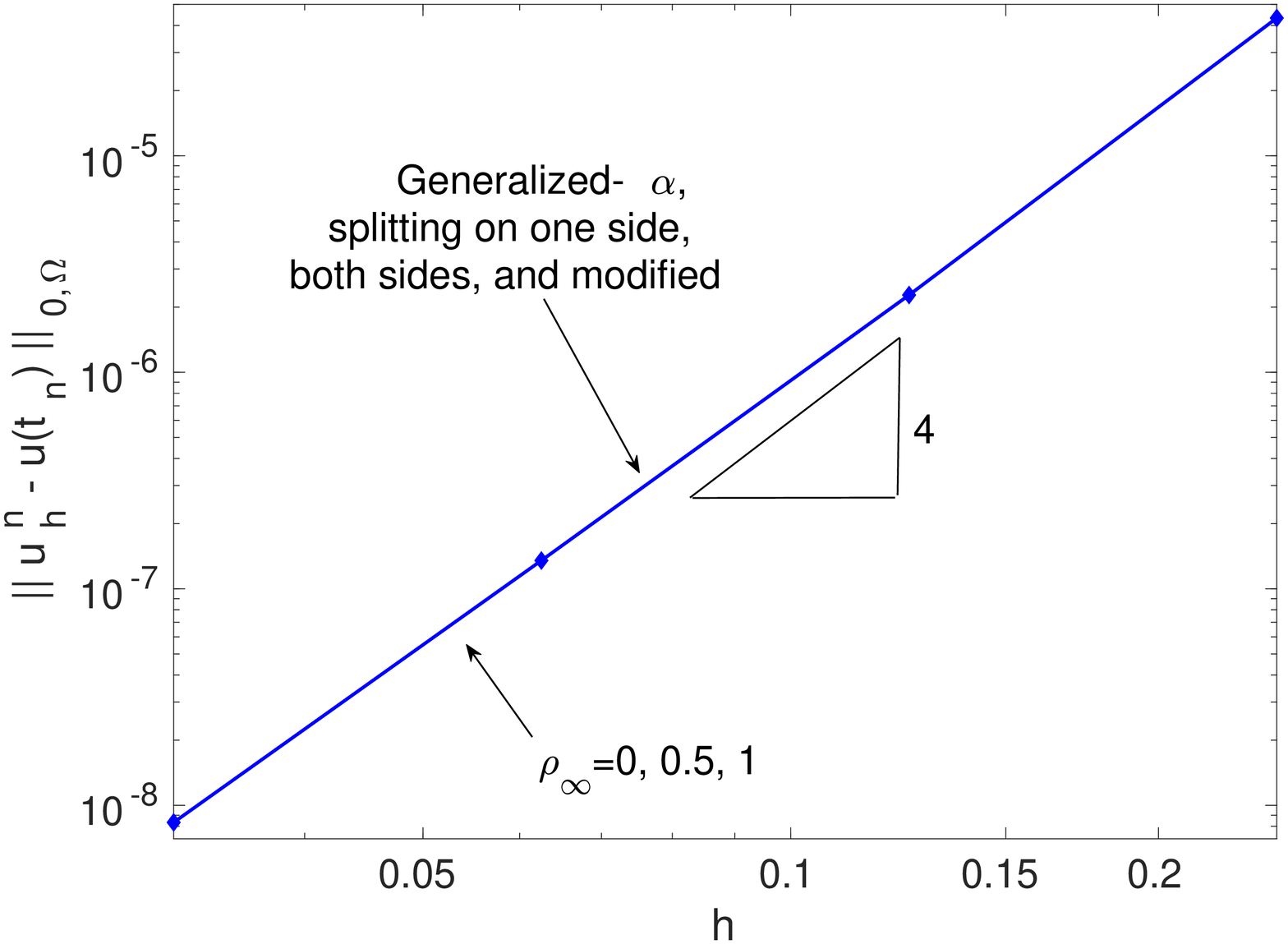} 
\hspace{0.3 cm}
\centering\includegraphics[width=6.5cm]{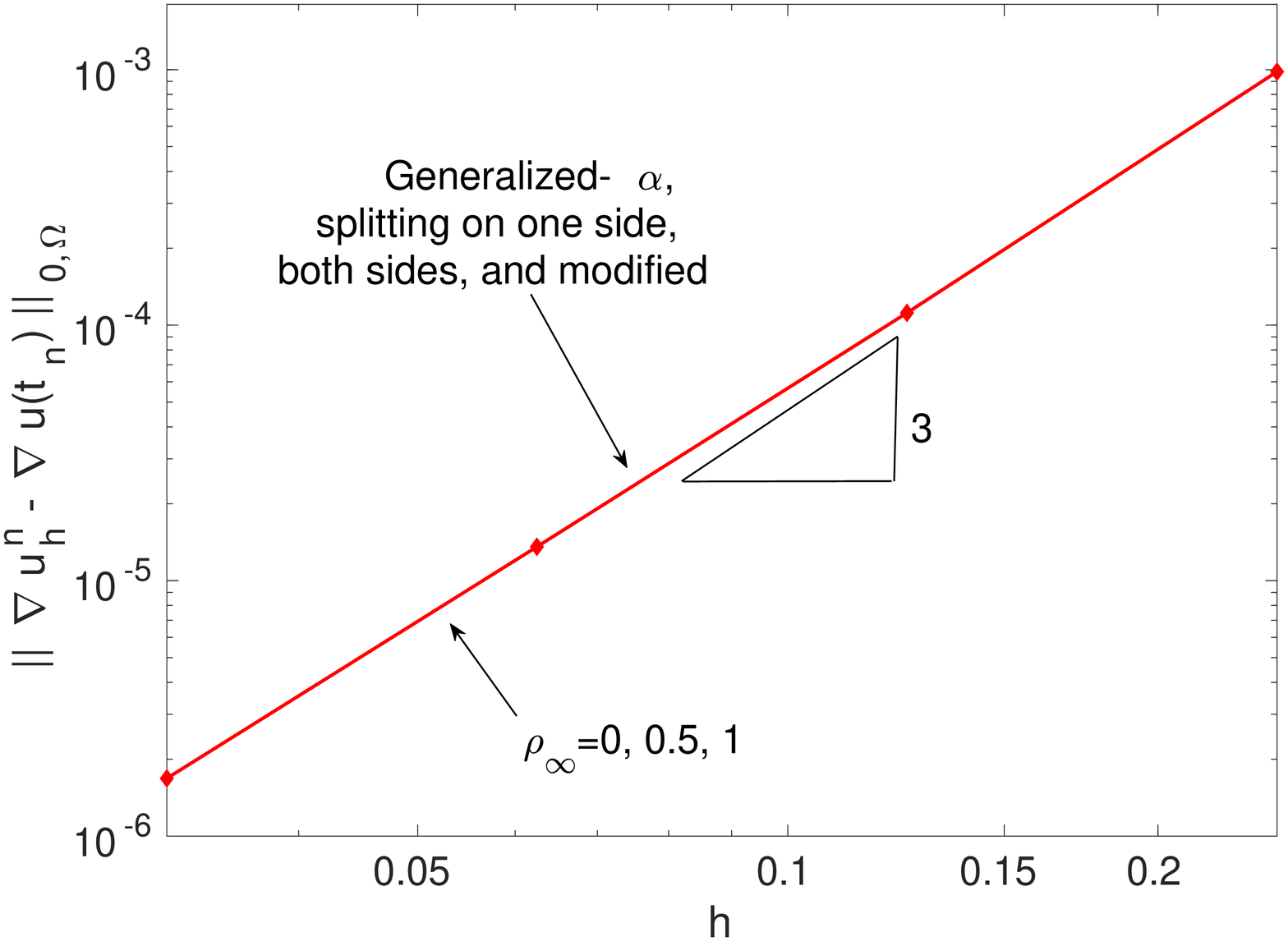} 

\caption{$L^2 $ norm of $u$ and $\nabla u$ errors when using  $C^2$ cubic isogeometric elements for space discretization with $\rho_\infty=0, 0.5, 1$. The final time $t_n=T=1$ and time step-size is $\tau = 10^{-4}$.}
\label{fig:p3h}
\end{figure}

Next, we verify the convergence rate of the spatial discretization error. The proposed splitting method can be applied to spatial discretizations using classical finite element analysis as well as isogeometric analysis. We consider the 2D test problem and fix the time step size to $\tau = 10^{-4}$. The final time for the simulation is $T =1$. Figure \ref{fig:p2h} shows the errors $\|u(T) - u^h(T)\|_{0,\Omega}$  and  $\| \nabla (u(T) -  u^h(T))\|_{0,\Omega}$ when using $C^0$ and $C^1$ quadratic elements. Figure \ref{fig:p3h} shows these errors when using $C^2$ cubic isogeometric elements. The generalized-$\alpha$ method and all the proposed splitting schemes result in optimal convergence rates for both finite element and isogeometric element discretizations. Clearly, the error in all cases converges with the corresponding optimal convergence rate.

\begin{figure}[!ht]
\subfigure[$\rho_\infty=0.5$, FEM $(p=2, C^0)$]{\centering\includegraphics[width=6.5cm]{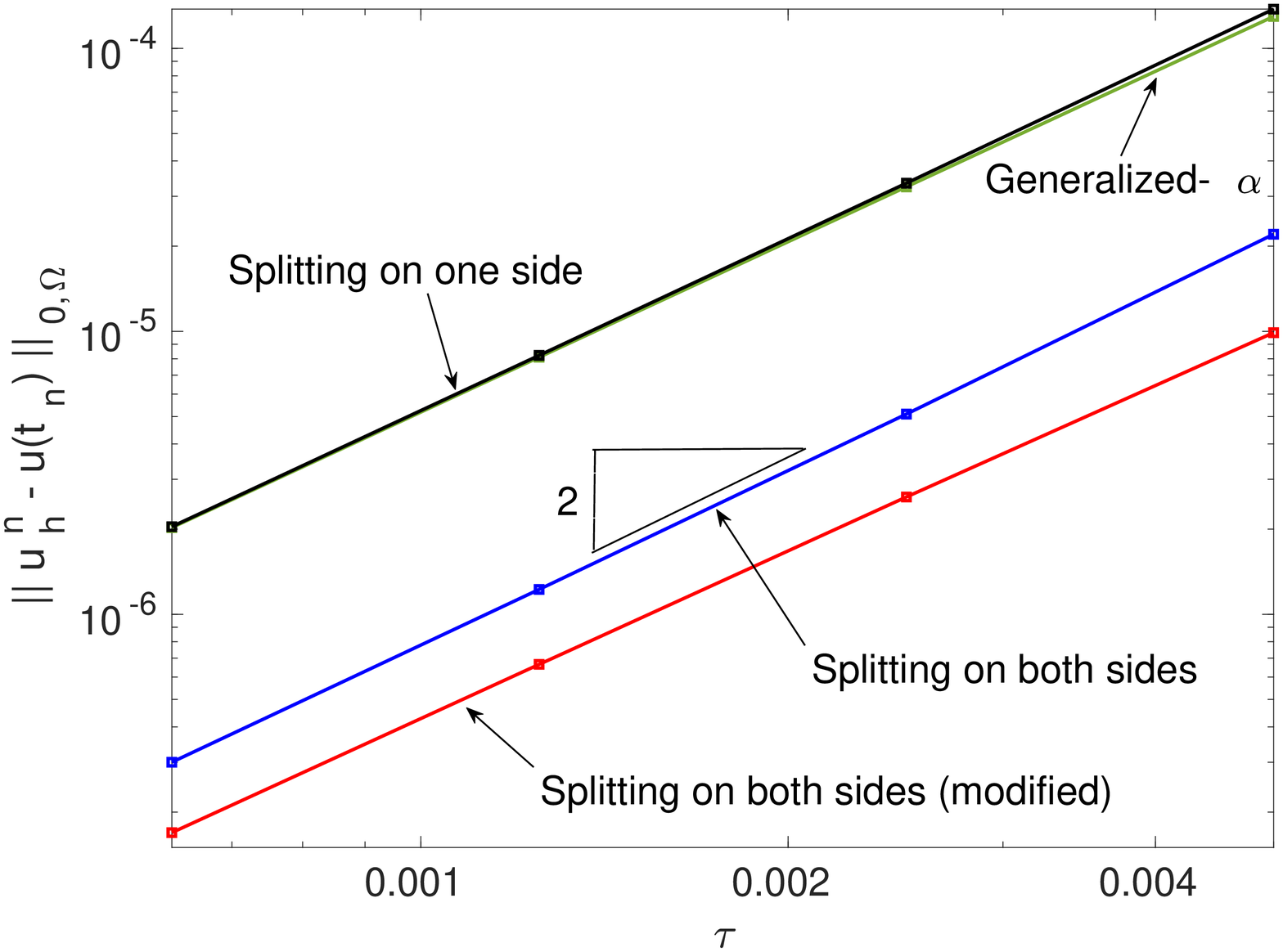} }
\subfigure[$\rho_\infty=0.5$, IGA $(p=2, C^1)$]{\centering\includegraphics[width=6.5cm]{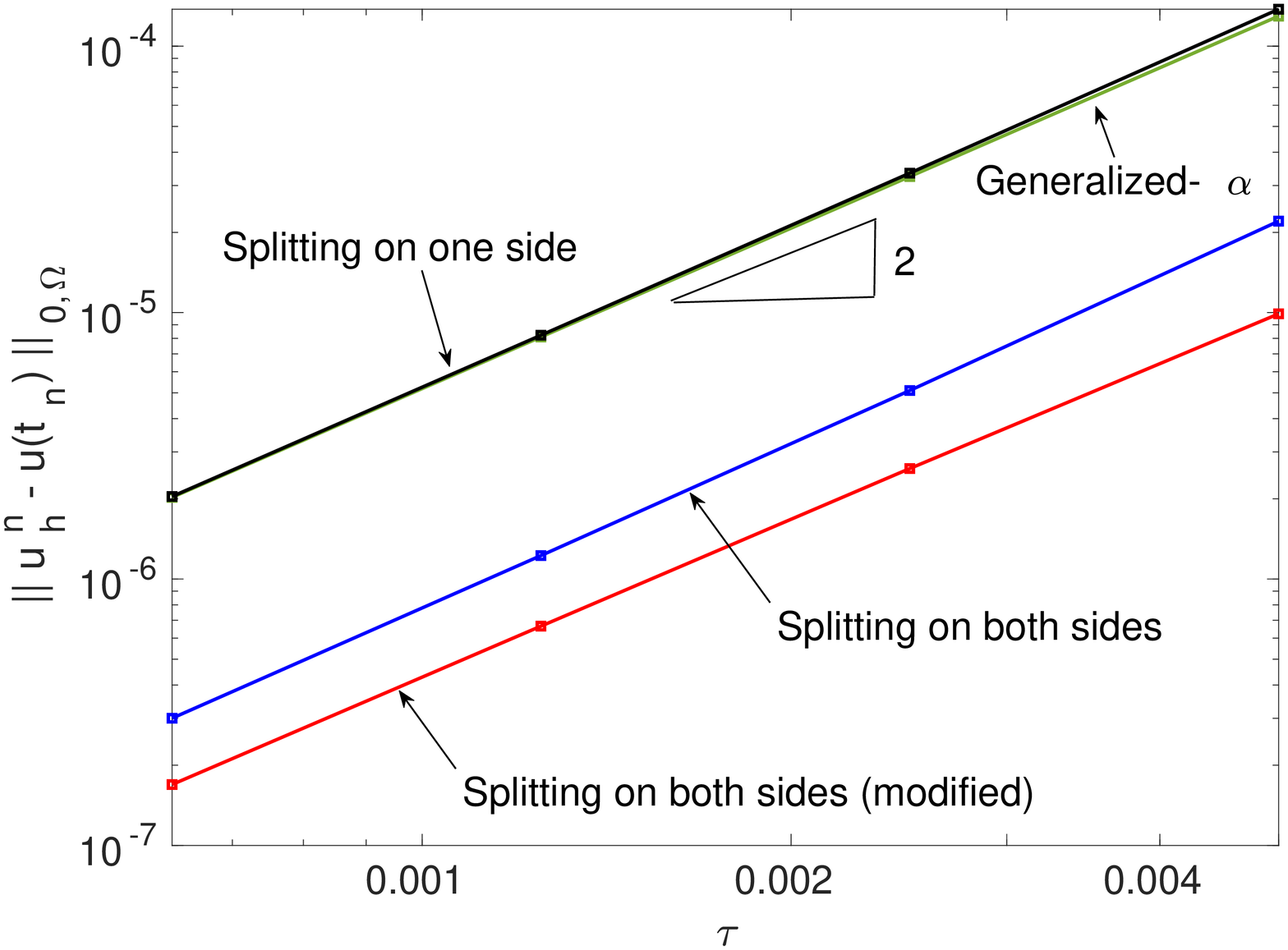} }
\subfigure[$\rho_\infty=1$, FEM $(p=2, C^0)$]{\centering\includegraphics[width=6.5cm]{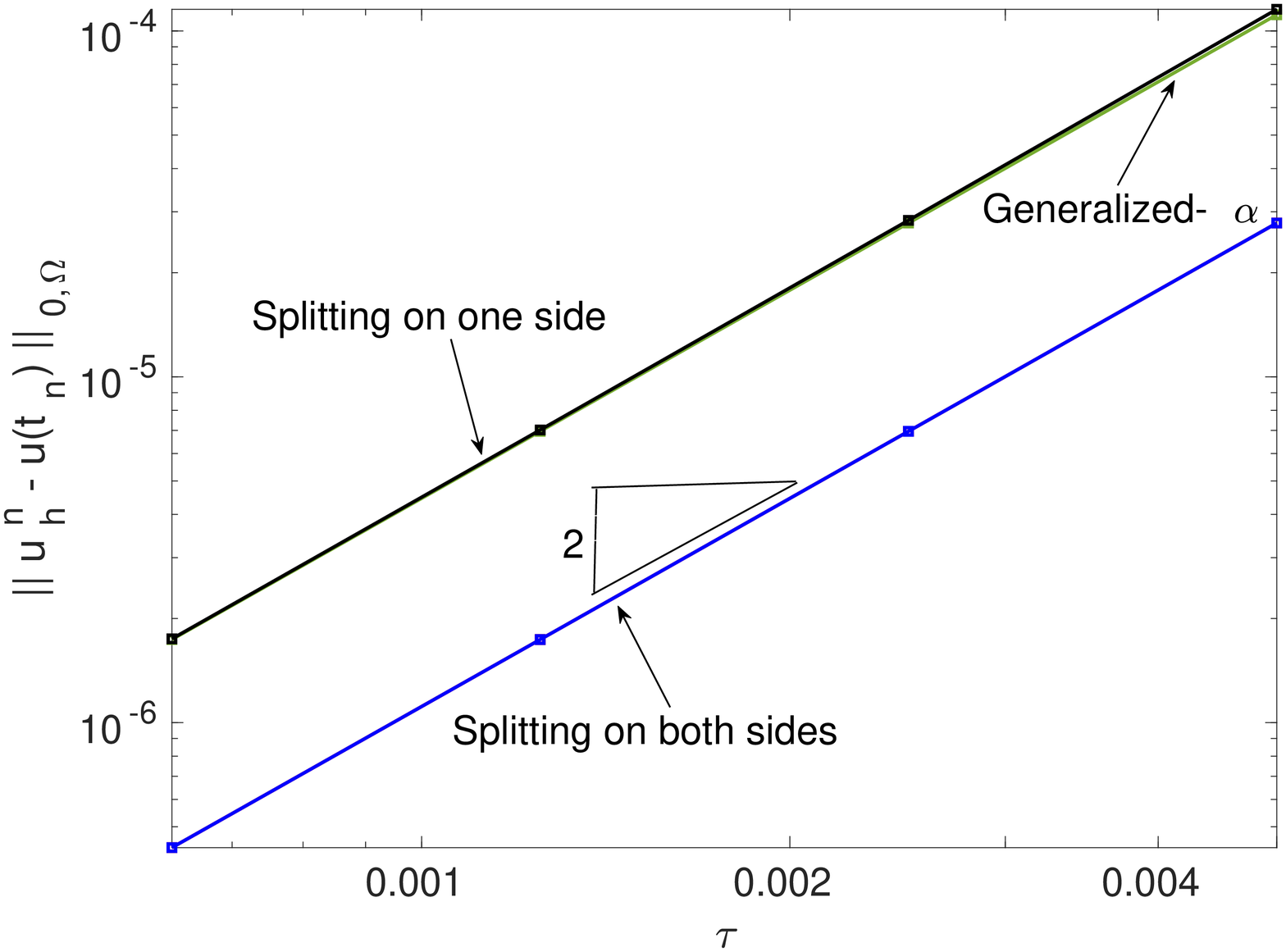}}
\hspace{0.12 cm}
\subfigure[$\rho_\infty=1$, IGA $(p=2, C^1)$]{\centering\includegraphics[width=6.5cm]{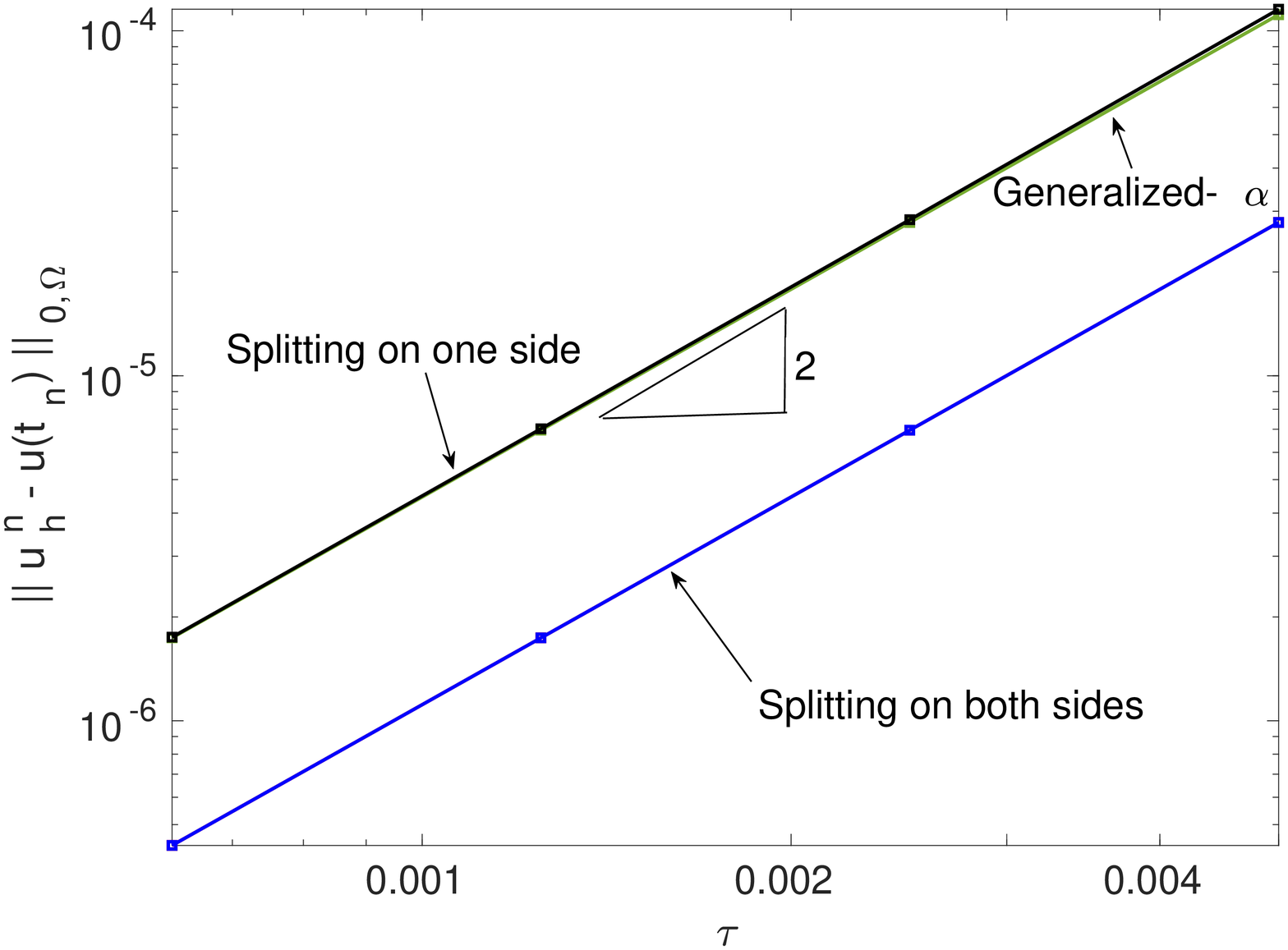}}
\caption{$L^2 $ norm error when using classical finite and isogeometric elements for space discretization with a fixed mesh size $h=1/64$ and $\rho_\infty=0, 0.5, 1$. The final time is $t_n = T =1.$}
\label{fig:p2t}
\end{figure}

\begin{figure}[!ht]
\subfigure[$\rho_\infty=0$]{\centering\includegraphics[width=6.5cm]{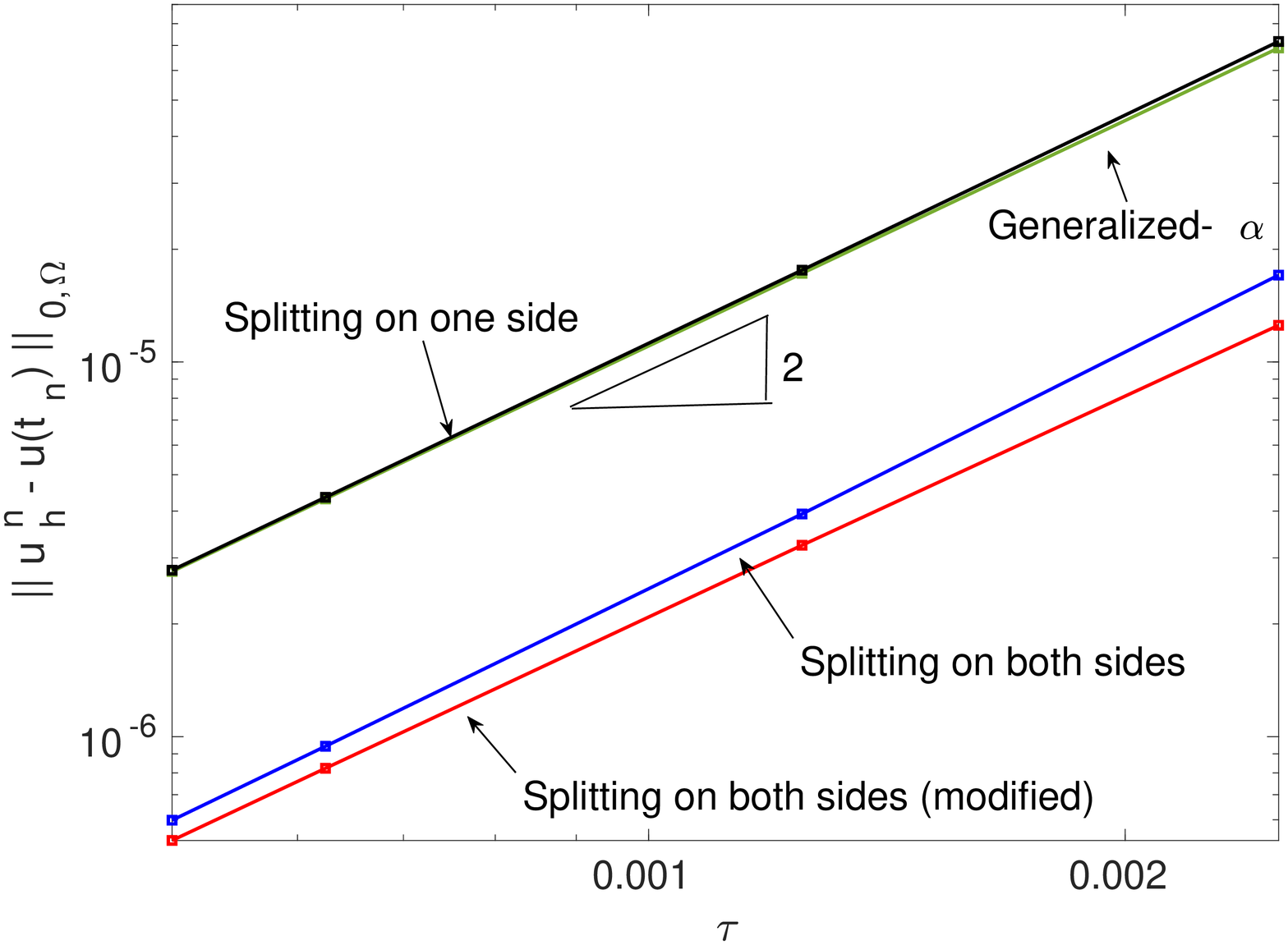} }
\subfigure[$\rho_\infty=1$]{\centering\includegraphics[width=6.5cm]{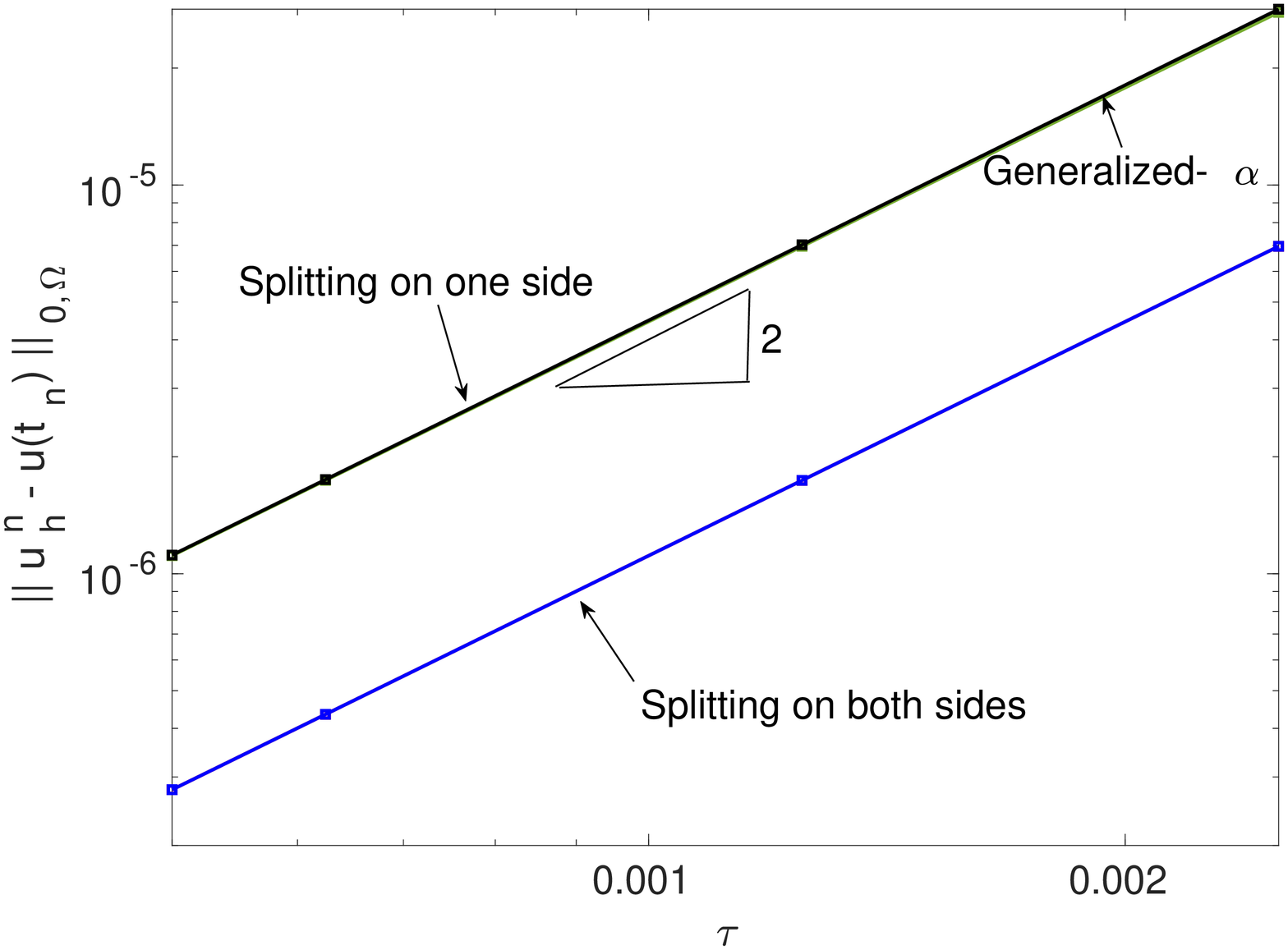} }
\caption{$L^2 $ norm error when using $C^2$ cubic isogeometric elements for space discretization with a fixed mesh size $h=1/64$ and $\rho_\infty=0, 0.5, 1$. The final time is $t_n = T =1$.}
\label{fig:p3t}
\end{figure}

Lastly, we study numerically the accuracy in time. We again consider the 2D test problem and use a mesh containing $100\times 100$ uniform elements. 
Figure \ref{fig:p2t} presents the $L^2$ norm of the  errors to show the convergence behaviour with respect to time step size $\tau$ when using quadratic finite elements and isogeometric elements. Figure \ref{fig:p3t} shows these errors when using $C^2$ cubic isogeometric elements. In  all cases, the errors converge quadratically thus verifying the formal estimates in \eqref{eq:errl2h1}.

\section{Concluding remarks}

We propose a splitting technique for parabolic equations using Cartesian product finite element or IGA basis  for spatial discretization, based on the generalized-$\alpha$ method for temporal discretization. The resulting splitting schemes are unconditionally stable and are second-order accurate in time while maintaining optimal convergence rates in space. The overall cost for solving the resulting algebraic system scales linearly with the number of degrees of freedom in the system,  and the computational cost is significantly reduced. These splitting schemes will be extended to hyperbolic equations and to other time integrators, and these developments will be reported in the near future.

\section*{Acknowledgement}
This publication was made possible in part by the CSIRO Professorial Chair in Computational Geoscience at Curtin University and the Deep Earth Imaging Enterprise Future Science Platforms of the Commonwealth Scientific Industrial Research Organization, CSIRO, of Australia. Additional support was provided by the European Union's Horizon 2020 Research and Innovation Program of the Marie Sk{\l}odowska-Curie grant agreement No. 777778, the Mega-grant of the Russian Federation Government (N 14.Y26.31.0013), the Institute for Geoscience Research (TIGeR), and the Curtin Institute for Computation. The J. Tinsley Oden Faculty Fellowship Research Program at the Institute for Computational Engineering and Sciences (ICES) of the University of Texas at Austin has partially supported the visits of VMC to ICES. The first and second authors also would like to acknowledge the contribution of an Australian Government Research Training Program Scholarship in supporting this research. Q. Deng thanks Professor Alexandre Ern (ENPC) for discussions on stability analysis. 
P. Minev acknowledges the support from a Discovery grant of  the Natural Sciences and Engineering Research Council of Canada, and makes an acknowledgment to the donors of the American Chemical Society Petroleum Research Fund for partial support of this research,  by a grant \# 55484-ND9.

%\section*{References}

\bibliographystyle{elsarticle-harv}\biboptions{square,sort,comma,numbers}
\bibliography{ref}

\end{document}